\documentclass[a4paper,fleqn,reqno]{cas-sc}
\usepackage[numbers, compress]{natbib}
\usepackage{amssymb, amsfonts}
\usepackage{soul}
\usepackage{subcaption}
\usepackage{placeins}

\def\tsc#1{\csdef{#1}{\textsc{\lowercase{#1}}\xspace}}
\tsc{WGM}
\tsc{QE}
\tsc{EP}
\tsc{PMS}
\tsc{BEC}
\tsc{DE}

\begin{document}
\let\WriteBookmarks\relax
\def\floatpagepagefraction{1}
\def\textpagefraction{.001}

\shorttitle{Broadband traveling wave generation using acoustic black holes}

\shortauthors{Omidi Soroor et~al.}

\title [mode = title]{Broadband traveling wave generation using acoustic black holes}                      


\author[1]{Amirhossein {Omidi Soroor}}[
type=editor,
 auid=000,bioid=1,
 orcid=0000-0002-9980-3924]

\cormark[1]


\ead{omidisoroor@tamu.edu}


\credit{Conceptualization, Methodology, Software, Validation, Formal analysis, Investigation, Data Curation, Writing - Original Draft, Writing - Review \& Editing, Visualization}

 \affiliation[1]{organization={J. Mike Walker '66 Department of Mechanical Engineering},
    addressline={Texas A\&M University}, 
    city={College Station},
    postcode={TX 77843}, 
    country={United States}}

\author[1,2]{Skriptyan N.H. Syuhri}

\credit{Conceptualization, Methodology, Investigation, Writing - Original Draft, Writing - Review \& Editing, Visualization}

\affiliation[2]{organization={Department of Mechanical Engineering},
    addressline={University of Jember}, 
    city={Jember},
    postcode={68121}, 
    country={Indonesia}}
  
\author[1]{Sourabh Sangle}


\credit{Software, Writing - Original Draft, Writing - Review \& Editing, Visualization}

\author[1]{Pablo A. Tarazaga}

\credit{Conceptualization, Resources, Writing - Review \& Editing, Supervision, Funding acquisition}

\cortext[cor1]{Corresponding author}

\begin{abstract}
This work studies the effectiveness of acoustic black holes to generate broadband non-reflective traveling waves using a single excitation source. This is inspired by similar observations in the basilar membrane of the mammalian inner ear. An aluminum beam is machined to introduce a gradual, asymmetric power-law taper at one of its ends. This tapered termination was then partially covered by viscoelastic tape to enhance the acoustic black hole effect in the system. This setup is then used to validate a model developed based on the Euler-Bernoulli beam theory. Following good agreement between the model and the experimental setup over a broad excitation frequency range ($1~kHz$ to $10~kHz$), the model is used to conduct a parametric study investigating the effects of different variables on the system's response. This study revealed the effectiveness of acoustic black holes in sustaining traveling waves over a broad frequency range. By optimizing parameters, such as the power-law order, one can significantly enhance this effect. This is especially noticeable towards the lower-frequency end. 
\end{abstract}

\maketitle
\section{Introduction}

Waves that continually progress from the point of incidence without any reflections are known as Traveling Waves (TWs). Reflections are due to impedance discontinuities along the path of the wave in a medium. As such, TWs are intuitively expected only to be observed in infinite media. However, there are examples in nature where TWs are observed in finite media in the form of microorganism motility \cite{taylor1952action,stone1996propulsion,vig2012swimming} and undulatory locomotion \cite{fish1982function,long1996importance,marvi2013snakes,chong2022coordinating}, to name a few. One of the interesting cases in which TWs naturally manifest is within the mammalian cochlea, specifically within a tiny biological structure called the basilar membrane (BM) \cite{von1928theorie,von1960experiments}. Note that although referred to as a membrane in the biological context, it mechanically acts as a beam. von Békésy observed that sound excites a wave in the BM that travels from its base towards the helicotrema near its apex without any observed reflections. He believed that the helicotrema was responsible for absorbing the reflections of incident waves \cite{von1970travelling}. These waves peak at locations along the BM that correspond to the frequency of the incoming sound \cite{zwislocki1980theory}. 

The progressive nature of these waves gives rise to a variety of engineering applications, examples of which include object transfer \cite{loh2000object,zhao2022requirements,rogers2023directed,mansouri2024toward}, structural health monitoring \cite{zhu2012piezoelectric,lugovtsova2019analysis}, flow control and drag reduction \cite{musgrave2019turbulent,olivett2021flow,zhang2024achieving,ogunka2025simulations}, robotic locomotion \cite{chen2007design,qi2020novel,zhu2024double,ji2024bionic}, submarine propulsion \cite{malladi2017investigation,musgrave2021electro,gupta2024utilizing,syuhri2024travelling,hess2024continuum}, pumping mechanisms \cite{hernandez2013design,rajendran2023novel}, and impact localization \cite{alajlouni2018impact,alajlouni2020passive,alajlouni2022maximum}. Generation and control of TWs in continuous systems have been extensively discussed in the literature. The majority of reported work employs the two-point excitation method to generate such waves. In this approach, a finite structure is excited at two points with harmonic loads of the same frequency, with a phase lag corresponding to the excitation frequency. This produces waves that travel from the location of one actuator to another \cite{malladi2015characterization}. This dynamic behavior is attributed to partial impedance matching at the points of excitation and reflection cancellation at the boundaries \cite{avirovik2016theoretical}. This method is studied and tested on a variety of mechanical elements \cite{anakok2022study,zhu2021design,cheng2025traveling,musgrave2021guidelines,rogers2024experimental,rogers2024estimation,phoenix2015traveling}. In contrast,  a less common methodology replaces one of the actuators with an intermediate impedance discontinuity (e.g., a spring and a damper) that is tuned to the excitation frequency to produce the same effect \cite{Motaharibidgoli_2023}. This approach, referred to as the single-point excitation method, has been investigated for strings \cite{blanchard2015damping,cheng2017separation}, beams \cite{cheng2019co,Motaharibidgoli_2023,gupta2024exploring}, and pressure waves \cite{xiao2017separation}. There are also other methods discussed in the literature, such as the active sink method \cite{tanaka1991active}, which involves two actuators: one pumping energy into the system and the other actively absorbing it through impedance matching, leading to virtually semi-infinite conditions \cite{gabai2008generating}. 

In the present work, inspiration is drawn from the BM's dynamic behavior. The main goal is to develop a model capable of mimicking von Békésy's observations in the cochlea to produce biomechanical models that accurately represent the BM and could improve non-invasive investigations of the auditory system. However, the applications of this basic study are not limited to BM modeling and also benefit the other applications discussed earlier. Accordingly, the single-point excitation method is adopted since the sound entering the ear canal is the sole source instigating the bending waves on the BM. However, so far, none of the developed single-excitation models have been able to passively produce broadband TWs without the need to adjust the impedance discontinuity to the input load's frequency. To achieve this robustness, Omidi Soroor and Tarazaga \cite{omidi2023investigation,soroor2025non} attempted to build the tonotopic nature of the BM into their model. They suspected tonotopy to be the primary cause of the broadband TW behavior on the BM. Thus, they investigated whether a non-uniform distribution of stiffness through material and geometric grading could improve the model's behavior in the telephony frequency range ($300~Hz-3.4~kHz$). However, it turns out that for any given system's configuration, as one sweeps through the frequency range, a node coincides with the intermediate spring-damper system, neutralizing it at the corresponding frequency such that no TWs form. The same goal is pursued in this work, though from a different perspective.

The effectiveness of Acoustic Black Holes (ABHs) in promoting wide-band TWs in the single-point excitation model is studied. ABHs are a relatively new development that serves as an effective passive solution for broadband vibration attenuation, particularly targeting lightweight structures \cite{gautier2020recent}. The concept was first introduced by Mironov \cite{mironov1988propagation} in 1988. He concluded that a sufficiently smooth power-law transition in the thickness of a finite plate to zero would cause the wave speed to vanish eventually; hence, it cannot reflect back at the boundary since it never reaches it. As this is a merely mathematical phenomenon and practically impossible, Krylov \cite{krylov2004new} suggested using a partial free (unconstrained) layer damping treatment, i.e., a thin viscoelastic material (VEM) layer bonded to the ABH section, which compensates for the truncated ABH termination and improves its damping performance. This was experimentally demonstrated by Krylov and Winward \cite{krylov2007experimental}, who suggested that using a relatively thin piece of damping layer would significantly improve the damping performance of the ABH, a favorable outcome for lightweight applications. From this point onward, interest in the so-called "ABH effect" in various engineering applications has grown significantly. These include but are not limited to broadband vibration attenuation \cite{bayod2011application,bowyer2014damping,zhang2024vibroacoustic}, broadband energy harvesting \cite{zhao2014broadband,ji2019enhancement, du2024semi}, metamaterials' performance enhancement \cite{bezanccon2024thin, zhang2025multi}, and artificial cochlea design \cite{foucaud2014artificial}.

The focus of this study is to inspect ABHs as a passive, frequency-independent solution for reflected wave cancellation, aiming to mimic the BM's vibratory behavior. As mentioned earlier, the approaches discussed in the literature are not applicable for this purpose, as wave absorbers must be actively tuned to the frequency of the incoming incident wave. In the following sections, the system's mathematical model will be developed, and the solution methodology will be explained. Next, the experimental procedure associated with this work will be discussed. Based on appropriate standards, the VEM will be characterized and then used as a partial free-layer damping agent for an aluminum beam with a machined ABH termination. This setup is used to validate the developed mathematical model. Finally, the validated model serves as the basis for a parametric study to reveal the effects of system variables—including excitation frequency, VEM properties, and ABH geometry—on the system's overall response. The results indicate that ABHs are a promising solution to generate broadband TWs in finite media passively.

\section{Problem Formulation and Solution Methodology}
\label{sec:Section2}

\begin{figure}
  \centering
  \includegraphics[width=.8\textwidth]{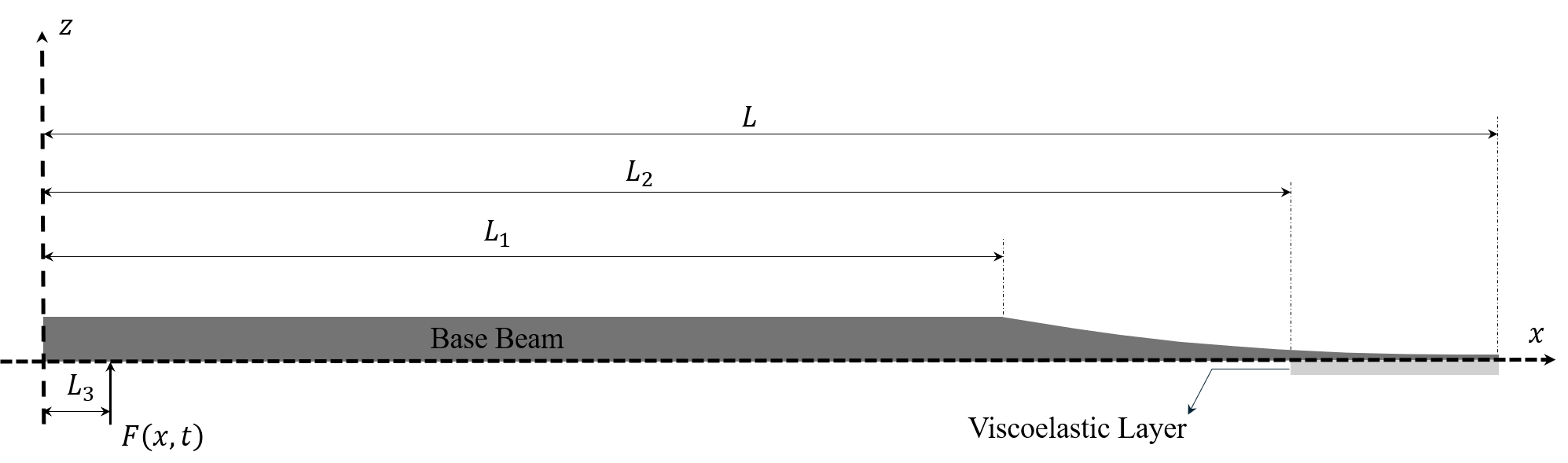}
  \caption{Beam with ABH termination and viscoelastic damping layer.}
  \label{fig:1}
\end{figure}

A free-free beam with an ABH termination coupled with a VEM damping layer is depicted in Fig.~\ref{fig:1}, where $L$ is the total length of the beam. $L-L_1$ and $L-L_2$ represent the ABH section and VEM layer lengths, respectively. $L_3$ is the location where the point load is exerted. It is assumed that this load is exerted at the midpoint of the piezoelectric transducer used in the experimental setup. $F(x,t)$ is the forcing term.

 Hamilton's principle, in the following form, will be used to develop the governing equation of the system as, 

\begin{equation}  
      \delta\int^{t_2}_{t_1}\bigg(T-U+W_{nc}\bigg)=0,
  \label{eq:1}
\end{equation}

\noindent where $T$, $U$, and $W_{nc}$ are the symbols for kinetic energy, strain energy, and virtual work due to non-conservative external forces, respectively. Here, $\delta$ denotes the variation operator, while $t_1$ and $t_2$ represent the initial and final instants of time. The Euler-Bernoulli beam assumption is adopted to model the system shown in Fig. \ref{fig:1}. The neutral axis of the composite beam, $\overline{z}$, at any given cross-section can be obtained as follows,

\begin{equation}
    \overline{z}=\frac{\sum_i{E_iA_i(x)z_i}}{\sum_iE_iA_i(x)},
    \label{eq:2}
\end{equation}

\noindent where $i$ refers to the subscripts $b$ and $v$ for the base and VEM damping layers, respectively. $z_i$ is layer $ i$'s centroid's location along the $z$-axis in the global system of coordinates shown in Fig. \ref{fig:1}. $E_i$ and $A_i$ represent the Young's modulus and cross-sectional area of layer $i$. Expanding Eq. \ref{eq:2} yields the neutral axis in the following form,

\begin{equation}
    \overline{z}=\frac{E_bh_b^2-E_vh_v^2}{2\big(E_bh_b+E_vh_v\big)},
    \label{eq:3}
\end{equation}

\noindent where, $h_i$ is the thickness of layer $i$. Finally, the equivalent bending stiffness for the composite beam, $D$, is calculated as follows,

\begin{equation}  
      D(x) =E_b\Big(I_b(x)+A_b(x)\big(z_b-\overline{z}\big)^2\Big)+E_v\Big(I_v(x)+A_v(x)\big(z_v-\overline{z}\big)^2\Big).
  \label{eq:4}
\end{equation}

In Eq. \ref{eq:4}, $I_i$ is the centroidal moment of inertia for each layer, i.e., $Bh_i^3/12$. $B$ is the beam's width. The thickness profile of the base beam is formulated as follows,

\begin{equation}
      h_b(x)=h_1\Big(H\big(x\big)-H\big(x-L_1\big)\Big)+\Bigg(h_2+\bigg(\frac{L-x}{L-L_1}\bigg)^m\big(h_1-h_2\big)\Bigg)\Big(H\big(x-L_1\big)-H\big(x-L\big)\Big).
  \label{eq:5}
\end{equation}

In Eq. \ref{eq:5}, $H(x)$ is the Heaviside function. Also, the VEM layer thickness across the full length of the beam is defined as follows,

\begin{equation}
      h_v(x)=h_3\bigg(H\big(x-L_2\big)-H\big(x-L\big)\bigg).
  \label{eq:6}
\end{equation}

In Eqs. \ref{eq:5} and \ref{eq:6}, $h_1$, $h_2$, and $h_3$ are values associated with the thickness of the uniform section of the base beam, the thinnest part of the ABH section, and the VEM layer, respectively. Additionally, $m$ is the power-law order. The strain energy of this system can then be calculated as follows,

\begin{equation}
      U =\frac{1}{2}\int_0^LD(x)w_{,xx}^2dx,
  \label{eq:7}
\end{equation}

\noindent where $w$ represents the beam's lateral displacement. Next, the kinetic energy expression is obtained as follows,

\begin{equation}
      T =\frac{1}{2}\int_{0}^L\mu(x)\dot w^2dx
      =\frac{B}{2}\int_{0}^L\big(h_b(x)\rho_b+h_v(x)\rho_v\big)\dot w^2dx.
  \label{eq:8}
\end{equation}

In Eq. \ref{eq:8}, $\mu(x)$ is the mass per unit length for the composite beam and $\rho_i$ is the layer $i$'s density. Finally, the variation in the virtual work is

\begin{equation}
      \delta W_{nc} =\int_0^LF\delta(x-L_3)\delta wdx.
  \label{eq:9}
\end{equation}

By inserting the expressions given in Eqs. \ref{eq:7}-\ref{eq:9} into Eq. \ref{eq:1} with some mathematical manipulation, the governing equation of the system is obtained as follows,

\begin{equation}
      \big(Dw_{,xx}\big)_{,xx}
      +\mu\ddot{w}=F\delta\big(x-L_3\big).
  \label{eq:10}
\end{equation}

The Galerkin method \cite{RAO} is implemented to discretize Eq.\ref{eq:10}. As such, a solution of the following form,

\begin{equation}
      w(x,t) =\sum_{l=0}^n\tau_l(t)\varphi_l(x), 
  \label{eq:11}
\end{equation}

\noindent is assumed, where $\varphi(x)$ and $\tau(t)$ represent the spatial trial functions and generalized coordinates. The Legendre polynomials, $P_n(x)$, are used as the trial functions in the following form,

\begin{equation}
  \varphi_l\big(\xi(x)\big)\;=\;
  \left\{
    \begin{aligned}
      &P_0(\xi) = 1 \\
      &P_1(\xi) = \xi \\
      &P_{n+1}(\xi) = \frac{1}{n+1}\bigl((2n+1)\xi P_n(\xi)-nP_{n-1}(\xi)\bigl),~l\geq2
    \end{aligned}
  \right.,
  \label{eq:12}
\end{equation}

\noindent where $\xi=2x/L-1$ and $l=n+1$. It is worth mentioning that by plugging the assumed solution given in Eq. \ref{eq:11} into the governing Eq. \ref{eq:10}, the discrete form of the latter is obtained as,

\begin{equation}
      \boldsymbol{M}\boldsymbol{\ddot{\tau}}+\boldsymbol{K}{\boldsymbol{\tau}}=\boldsymbol{f},
  \label{eq:13}
\end{equation}

\noindent where $\boldsymbol{M}$ and $\boldsymbol{K}$ are the mass and complex stiffness matrices and $\boldsymbol{f}$ is the forcing vector. Assuming that the point load is harmonic of the form $F=F_0e^{i\omega t}$, the generalized coordinates take a similar harmonic form, that is, $\boldsymbol{\tau}=\boldsymbol{\tau_0}e^{i\omega t}$. Therefore, the forced response of the system is obtained by using,

\begin{equation}
      \boldsymbol{\tau_0}=\Big[-\omega^2\boldsymbol{M}+\boldsymbol{K}\Big]^{-1}\boldsymbol{f_0},
  \label{eq:14}
\end{equation}

\noindent where $\boldsymbol{f_0}$ is the forcing amplitude vector. The matrices in Eq. \ref{eq:14} are derived using the weak form of the governing equation and are as follows,

\begin{equation}
\begin{aligned}
     &\boldsymbol{M} = \int_0^L\mu(x)\varphi_l\varphi_jdx\\[2pt]
     &\boldsymbol{K} =\int_0^LD(x)\varphi''_l\varphi''_jdx\\[2pt]
     &\boldsymbol{f_0}=F_0\varphi_j(L_3)
\end{aligned} .
  \label{eq:15}
\end{equation}

A code based on the discussed formulation has been developed and validated using the experimental investigation outcomes presented in the subsequent sections.

\section{Experimental Procedure}

The experimental procedure used in this study is presented below. In the first step, the VEM to be used as the free-layer damping treatment was characterized using ASTM E756-05 \cite{E756}, which is designated for the vibration-damping characterization of materials. Next, the modal properties of an Aluminum beam with a machined ABH section partially covered by VEM tape were acquired for model verification purposes. Finally, the dynamic response of the system due to harmonic loads of different frequencies within the valid range was obtained and compared against the model.

\begin{table}
\centering
\caption{Geometric properties of the substrate in the Oberst beam configuration}
\label{tab:oberst_properties}
\begin{tabular}{@{}ll@{}}
\toprule
\textbf{Parameter} & \textbf{Value ($mm$)} \\
\midrule
Root section length & $30.48$ \\
Root section width & $25.40$  \\
Root section thickness & \(9.525\) \\
Beam section length & $170.0$ \\
Beam section width & $25.40$  \\
Beam section thickness & \(1.588\) \\
\bottomrule
\end{tabular}
\end{table}
\subsection{Viscoelastic Material Characterization}
\label{subsec:visco_mat_char}

Complex moduli are used to describe the elastic properties of VEMs. In this study, the Young's modulus is of interest as the shear modulus does not appear in the developed unconstrained damping formulation. The Young's modulus is mathematically idealized as follows,

\begin{equation}
 E_v =  E_{v_s} + iE_{v_l},
  \label{eq:16}
\end{equation}

where, $E_{v_s}$ and $E_{v_l}$ are the frequency-dependent storage and loss moduli, respectively. Inserting the loss factor, $\eta=\frac{E_{v_l}}{E_{v_s}}$, into Eq. \ref{eq:16} yields,

\begin{equation}
 E_v =  E_{v_s}(1+i\eta).
  \label{eq:17}
\end{equation}

\begin{figure}
\begin{subfigure}[b]{0.45\textwidth}
  \centering
  \includegraphics[width=\linewidth]{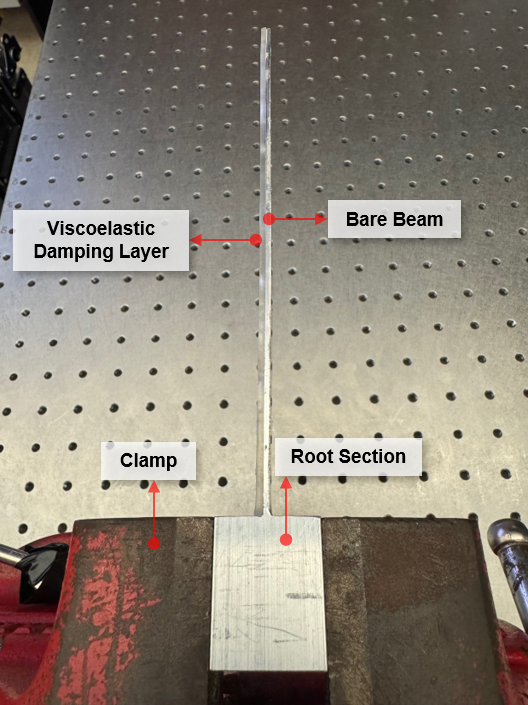}
  \subcaption{}\label{fig:damp_a}
\end{subfigure}
\hfill
\begin{subfigure}[b]{0.45\textwidth}
  \centering
  \includegraphics[width=\linewidth]{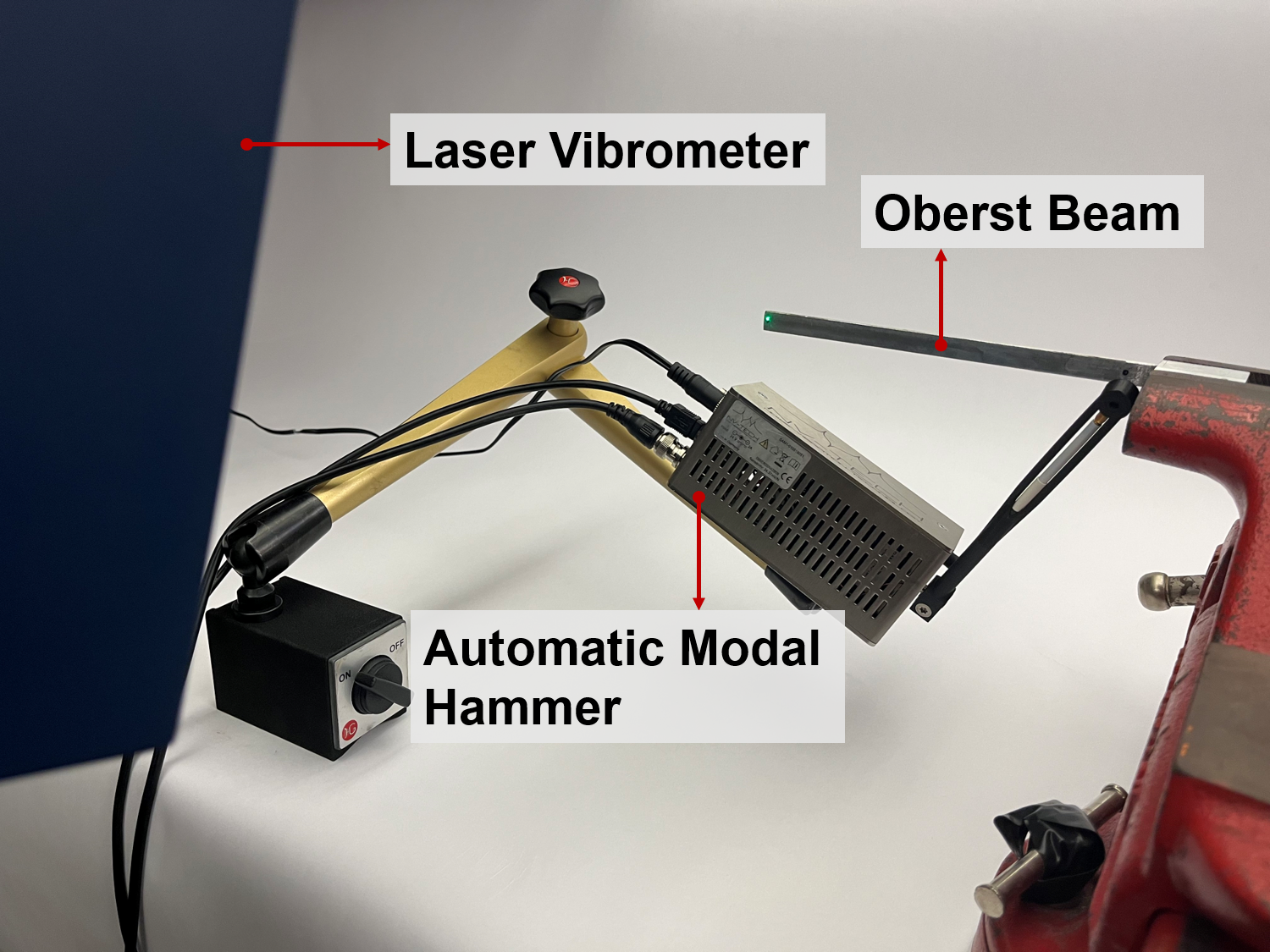}
  \subcaption{}\label{fig:damp_b}
\end{subfigure}
\caption{Experimental setup for viscoelastic damping characterization: (a) Oberst beam specimen with VEM layer, root section, and clamp; (b) dynamic testing with laser vibrometer and automatic modal hammer.}\label{fig:dampchar2}
\end{figure}

In Eq. \ref{eq:17}, the loss factor directly correlates to the damping performance of a VEM. The VEM used in this study was characterized following ASTM E756-05, using an Oberst beam configuration comprising a 6061-T4 aluminum substrate bonded with a damping layer as shown in Fig.~\ref{fig:dampchar2}(a). The specimen consisted of a root and a beam section, fabricated according to the standard procedure. Geometric properties of the substrate are summarized in Table~\ref{tab:oberst_properties}. The damping layer had a density of $1041.2~kg/m^3$ and a thickness of $1.7~mm$, and was trimmed to match the beam section dimensions after bonding to one side of the substrate.

Fig.~\ref{fig:dampchar2}(b) shows the experimental setup used to characterize the damping properties of the VEM. The setup consisted of the aforementioned manufactured beam to mimic a fixed-free boundary condition, an automatic modal hammer\footnote[1]{Scalable Automatic Modal Hammer type SAM1} for excitation, and a Scanning Laser Doppler Vibrometer (SLDV) system\footnote[2] {Polytec PSV QTec} for response measurement. The SLDV system incorporates a scanning head, signal generation and processing units with multiple I/O channels, and dedicated software for control and data acquisition.

The beam was excited at a location close to the root section, with a consistent impact position maintained throughout the tests. Although this method involves mechanical contact, it provides a repeatable excitation force, offering a practical compromise in place of the ideal non-contact excitation recommended by ASTM E756-05. The vibration responses were measured at the beam's free end. For each test, ten consecutive impacts were conducted, and the Polytec software automatically averaged the corresponding Frequency Response Functions (FRFs). Furthermore, all tests were performed under ambient laboratory conditions without thermal regulation.

\begin{figure}
    \centering
    \includegraphics[width=.8\textwidth]{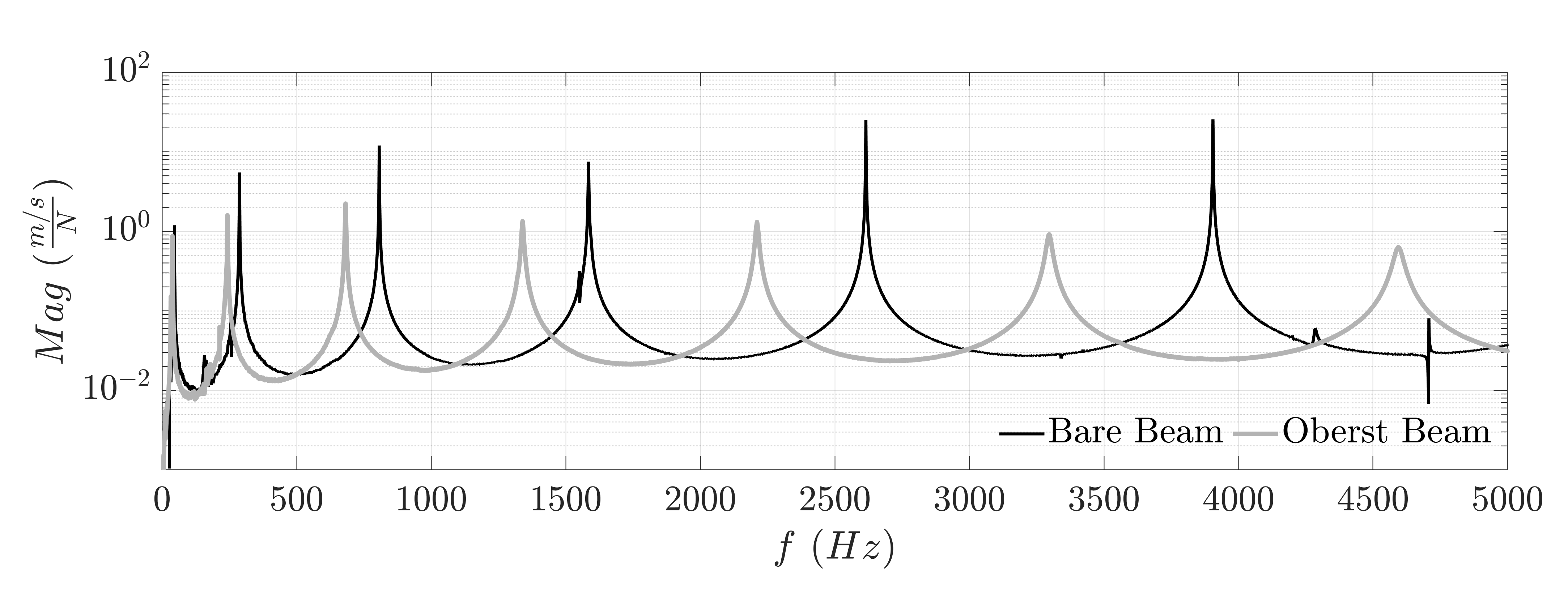}
    \caption{FRFs of the bare and Oberst beams measured at the free end, showing the effect of the damping layer on the resonance amplitude, frequency shift, and bandwidth.}
    \label{fig:frf_damping}
\end{figure}

Fig.~\ref{fig:frf_damping} presents the FRFs of the bare beam and Oberst beam configuration measured at the free end. A total of six resonance peaks are identified in the bare beam up to 5~kHz, while the Oberst beam shows seven peaks over the same frequency range. The addition of the viscoelastic damping layer leads to a noticeable reduction in peak amplitudes, particularly at higher frequencies, along with wider resonance peaks that indicate enhanced energy dissipation. The downward shift in resonant frequencies is most likely attributed to the added mass introduced by the damping layer. These trends demonstrate the effectiveness of the damping in attenuating vibration across a broad frequency range.

\begin{figure}
    \centering
    \includegraphics[width=.85\textwidth]{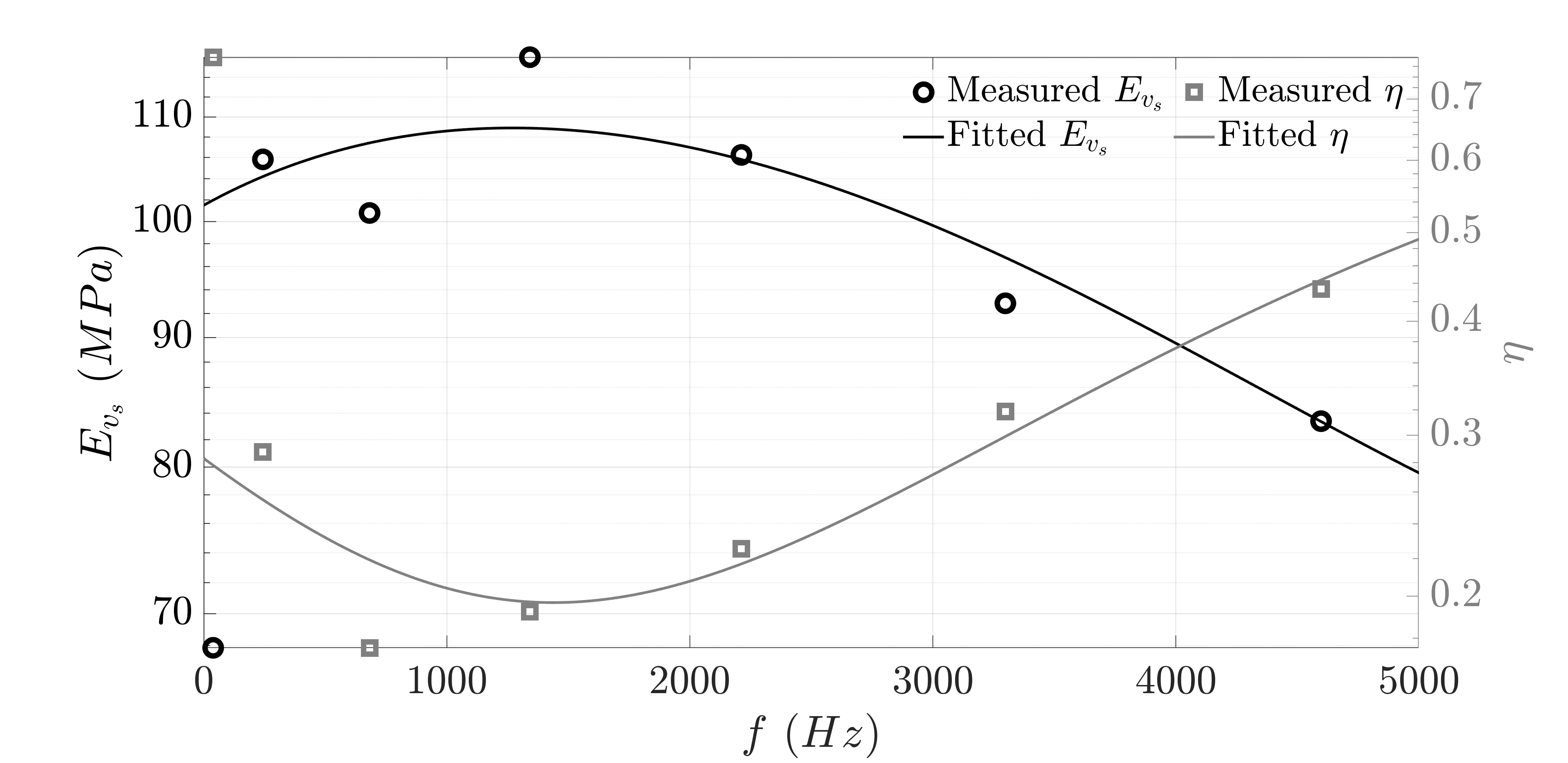}
    \caption{Frequency-dependent storage modulus and loss factor of the damping material, extracted from FRFs using ASTM E756-05.}
    \label{fig:Eandeta}
\end{figure}

The storage modulus and loss factor of the damping material were obtained from the FRFs of the bare and Oberst beams using the methodology described in the ASTM E756-05 standard. As shown in Fig.~\ref{fig:Eandeta}, the storage modulus decreases slightly, while the loss factor decreases with increasing frequency up to $5~kHz$, which is the recommended range in the used standard. This frequency-dependent reduction in dynamic stiffness is consistent with previously reported results for viscoelastic polymers. For example, \citet{Brun_2022} observed a similar decrease in modulus, from 75.6~GPa to 43.1~GPa, using a curve-fitting method based on transmissibility data. This behavior was consistent across multiple test cases and is attributed to internal energy dissipation mechanisms within the material. The decreasing trend observed in our study may, therefore, reflect similar mechanisms inherent to the composition of the viscoelastic tape used. 

Furthermore, in this study, the average values across the frequency range are used in the model, with a mean storage modulus of $96.16~MPa$ and a mean loss factor of $0.34$. This approach yields higher errors at lower frequencies \cite{soroor2021effect}, but achieves acceptable accuracy when higher frequencies are more dominant in the system response, which is the case in this work. As such, these average values are not expected to change as the frequency further increases due to the plateau in viscoelastic properties at higher frequencies that were not included in the material characterization \cite{denis2015measurement}. This is further confirmed in sections \ref{subsec:modal} and \ref{sub:MV}, where experimentally extracted results are compared to those obtained using the developed model. The error generally decreases as frequency increases.

\subsection{Test Specimen: Euler-Bernoulli Beam with ABH}
\label{subsec:ExpEB-ABH}

\begin{figure}
    \centering
    \includegraphics[width=0.95\textwidth]{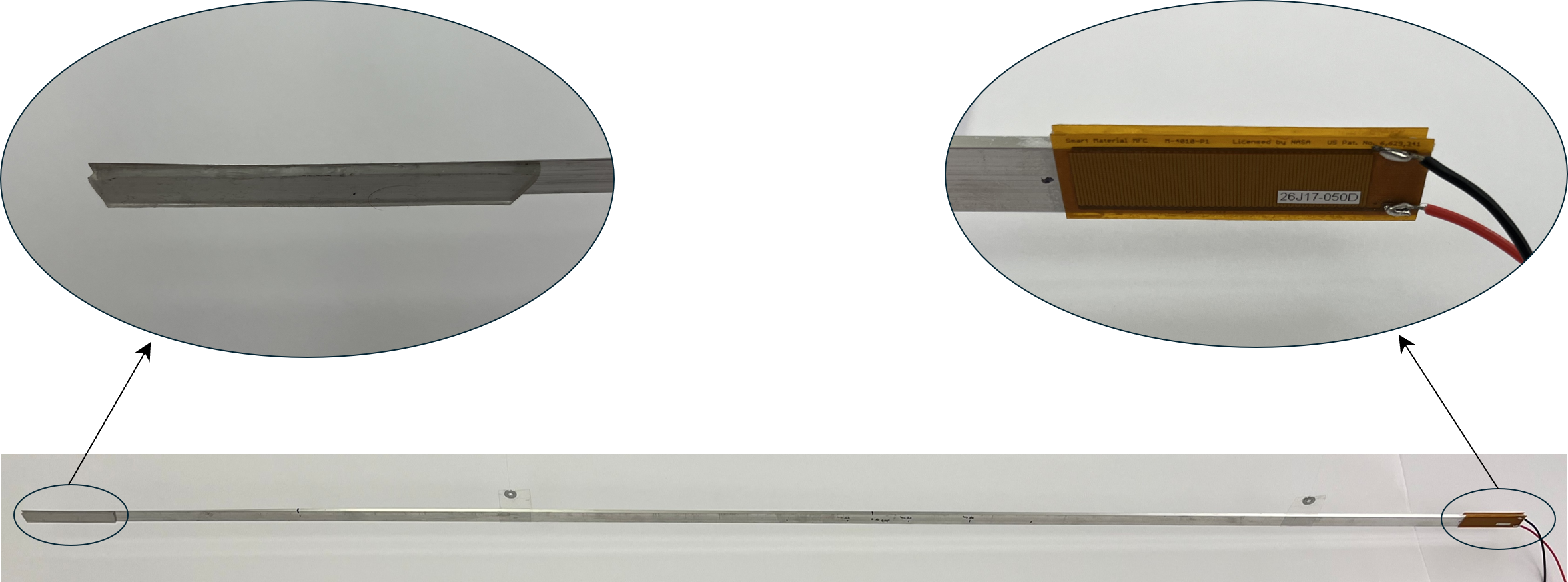}
    \caption{Free-free beam with an ABH, integrated with piezoelectric patches and a viscoelastic damping layer.}
    \label{fig:EBABHExp}
\end{figure}

Fig.~\ref{fig:EBABHExp} shows the beam with ABH in a free-free configuration. To accurately mimic this boundary condition, the beam was suspended using two thin transparent strings positioned at separate locations along the span. The beam was machined from 6063-T6 aluminum and consists of two sections: a uniform region and a tapered ABH region, with a total length of $1.22~m$. The ABH region adopts a power-law thickness profile with exponent $m = 3$, gradually tapering to a finite minimum thickness that ensures manufacturability. The values of $m$ and the tapered length $L-L_2$ were selected following the approach outlined in~\cite{hook2019}, where these parameters were shown to minimize wave reflection. To enhance damping performance, the previously characterized VEM was applied along roughly one-third of the ABH region, starting from the free tip. Moreover, a Smart Material M-4010-P1 Macro-Fiber Composite (MFC) piezoelectric patch was used near the free end of the uniform section to serve as an actuator for harmonic excitation. A complete summary of the geometric and material properties is provided in Table~\ref{tab:matpropEBABH}.

\begin{table}[width=.7\linewidth,cols=3,pos=h]
\caption{Summary of geometric and material properties of the beam with ABH.}\label{tab:matpropEBABH}
\begin{tabular*}{\tblwidth}{@{} LLL@{} }
\toprule
\textbf{Parameter} & \textbf{Value} & \textbf{Unit} \\
\midrule
\textbf{Beam (Aluminium 6063-T6)} & & \\
Density \( \rho \) & $2700$ & $kg/m^3$ \\
Modulus of elasticity \( E \) & $68.9$ & $GPa$ \\
Total length \( L \) & $1.22$ & $m$ \\
Beam width \( B \) & $1.27$ & $cm$ \\
Uniform section length \( L_1 \) & $1$ & $m$ \\
ABH section length \( L-L_1 \) & $22$ & $cm$ \\
Uniform thickness \( h_1 \) & $3$ & $mm$ \\
Minimum thickness \( h_2 \) & \(0.2\) & $mm$ \\
Thickness exponent \( m \) & $3$ & -- \\
\addlinespace
\textbf{Damping Layer (Double-sided adhesive tape)} & & \\
Thickness $h_3$ & \(1.9\) & $mm$ \\
Length $L-L_2$ & $8.2$ & $cm$ \\
\addlinespace
\textbf{Piezoelectric Actuator (MFC M-4010-P1)} & & \\
Length & $5$ & $cm$ \\
Width & $1.6$ & $cm$ \\
Thickness & \(0.3\) & $mm$ \\
\bottomrule
\end{tabular*}
\end{table}

\subsection{Experimental Modal Analysis}
\label{subsec:modal}

\begin{figure}
    \centering
    \includegraphics[width=0.6\textwidth]{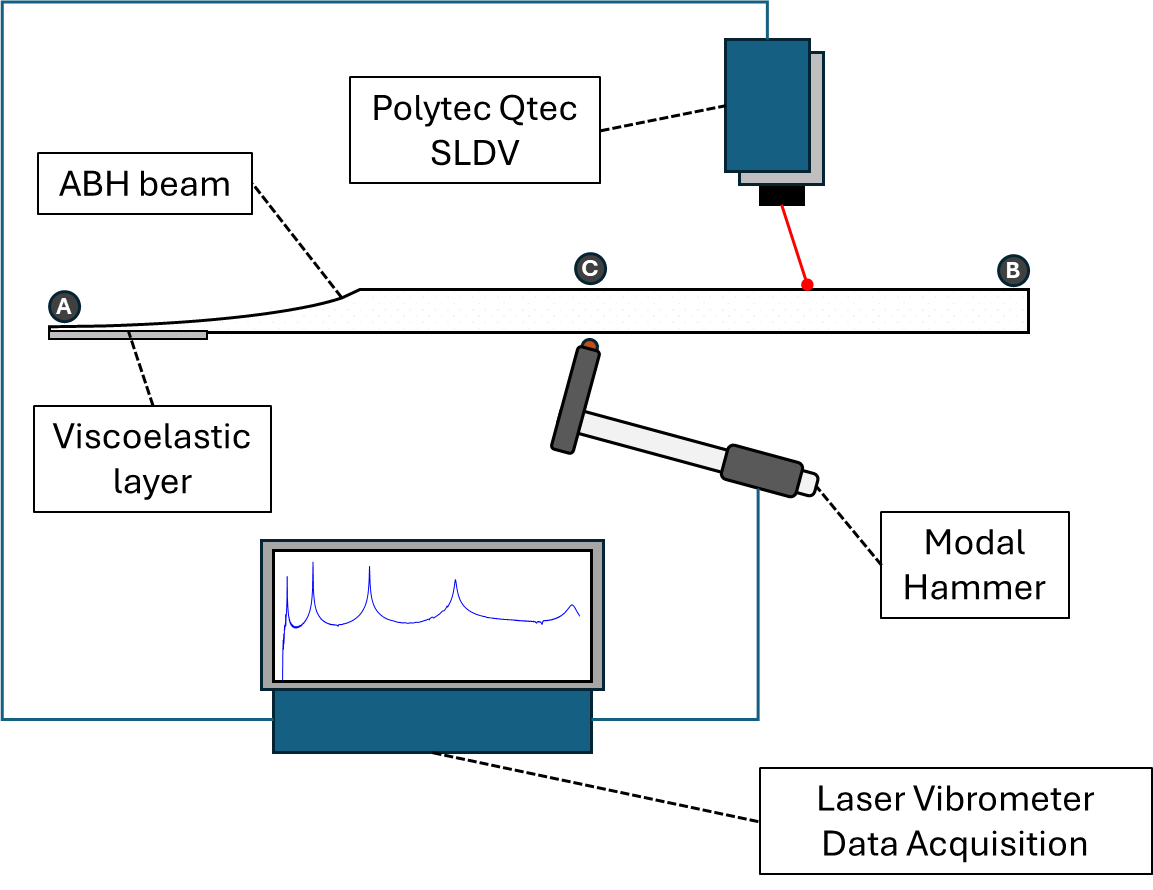}
    \caption{Schematics for experimental setup of the ABH beam with free-free configuration with points A, B, and C marked}
    \label{fig:exp_schematics}
\end{figure}

Fig. \ref{fig:exp_schematics} presents the schematics of the experimental setup for the free-free beam with an ABH. The same measurement equipment described in Section \ref{subsec:visco_mat_char} was employed. The excitation was applied at point C, as indicated in the figure, and maintained there throughout the test. The SLDV scan was performed along the beam surface, from point A to point B. Point A is located at the thin tip of the ABH section, and point B is just beyond the bonded piezoelectric patch. A total of 71 measurement points were performed to adequately capture the spatial resolution for this case. At each point, three individual impacts were conducted and subsequently averaged using the Polytec software to reduce the measurement noise.

\begin{figure}
    \centering
    \includegraphics[width=0.85\textwidth]{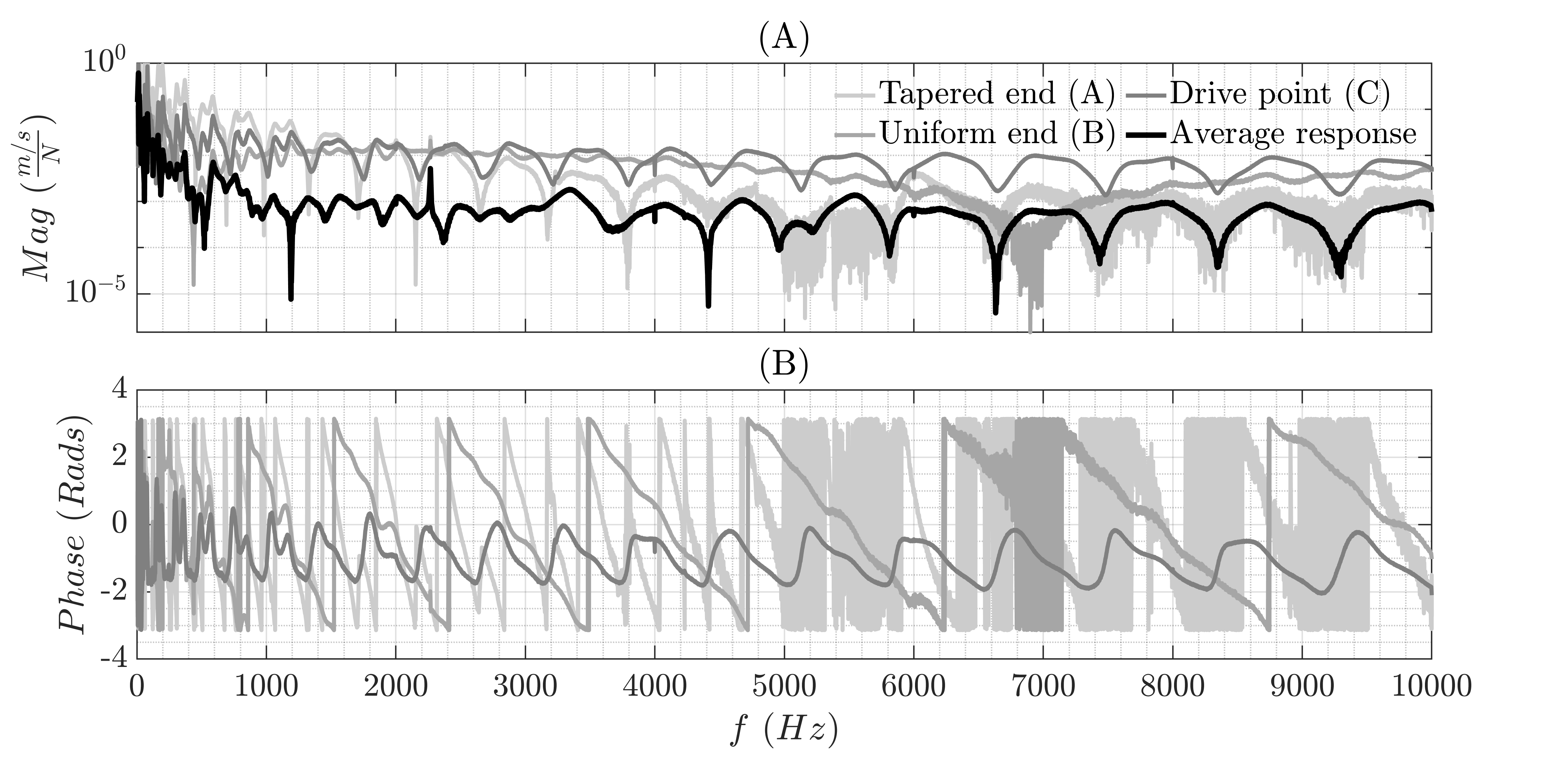}
    \caption{(A) Magnitude and (B) Phase for frequency response functions recorded at points A, B, and C.}
    \label{fig:frfs_3pts+avg}
\end{figure}

Figs. \ref{fig:frfs_3pts+avg}(A) and (B) present the magnitude and phase components of the measured FRF at three locations along the beam: the tapered end A, the uniform end B, and the excitation point C. The thick black line in the magnitude plot indicates the average FRF across the measurement points (excluding very noisy data). As observed, the system exhibits significant damping, which is evident from the broadening of the resonant peaks with increasing frequency. This damping effect reduces the clarity of peak locations, making it difficult to accurately extract natural frequencies using conventional peak-picking methods. To address this challenge, the Vector Fitting (VF) algorithm \cite{gustavsen_rational_1999,gustavsen_improving_2006,deschrijver_macromodeling_2008} was employed to build a reduced-order model and estimate the modal parameters. VF was selected due to its accuracy in approximating high-frequency response and its ability to capture the modal parameters via conjugate pole-pair estimation \cite{sangle2024pole}. Interested readers can refer to Ref. \cite{sangle2024pole} for more details about the advantages of VF over other techniques and its implementation.

\begin{figure} 
    \centering
    \includegraphics[width=0.85\textwidth]{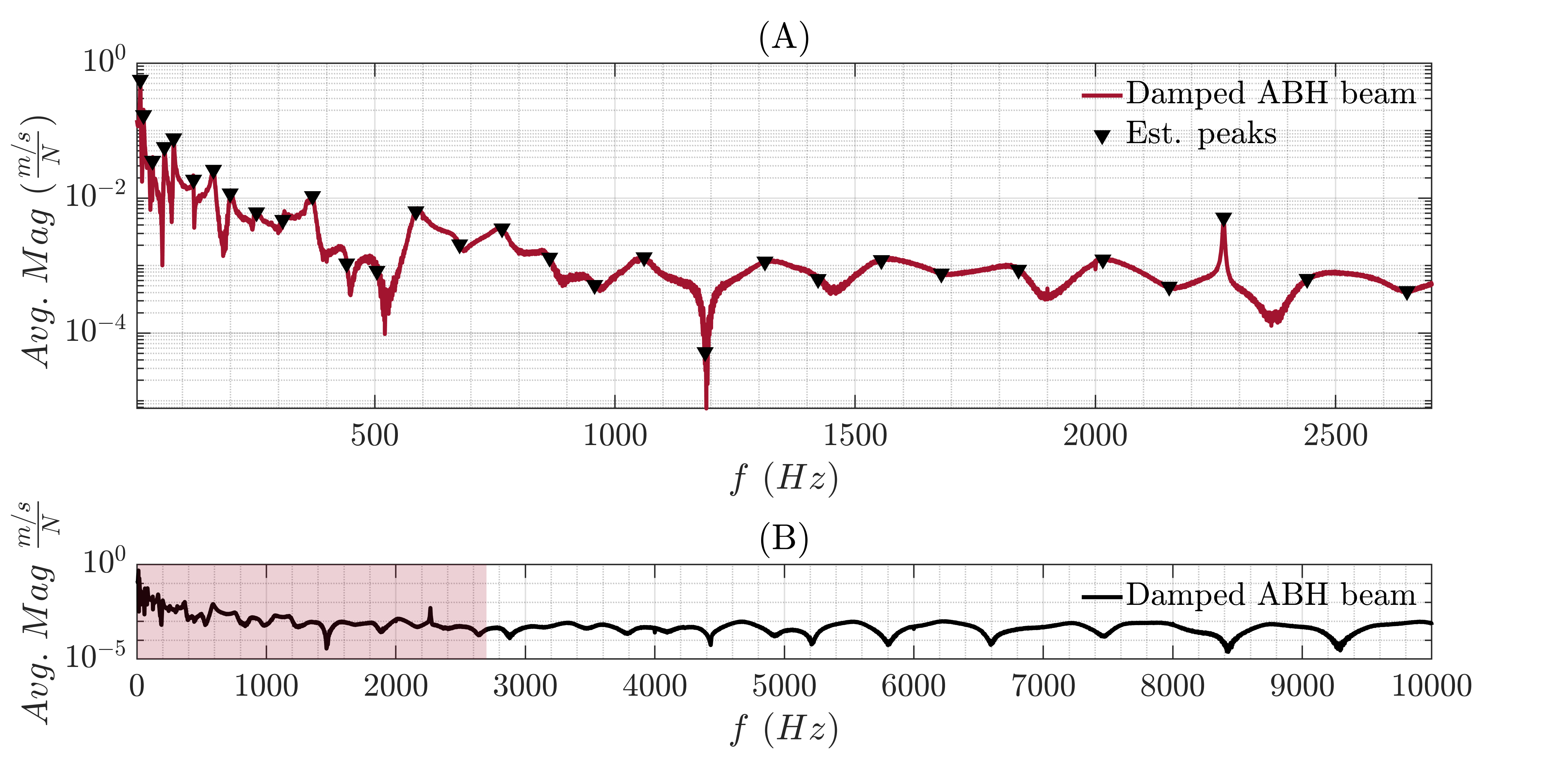}
    \caption{Average frequency response for ABH beam (A) over selected range of 0 - 2700 Hz with estimated peaks and (B) over entire range of 0 - 10 kHz}
    \label{fig:frfs_vf_peaks}
\end{figure}

In implementing VF, the frequency response functions were divided into smaller bands, and VF was employed to estimate the peaks (modes) over each. It is worth noting that the complex-valued frequency response functions were processed simultaneously for all points. Following the procedure described in Ref. \cite{sangle2024pole}, a set of frequencies was identified such that corresponding poles are stable and non-spurious. For this study, the analysis focused on capturing the first 30 bending modes using the average FRF shown in Fig. \ref{fig:frfs_vf_peaks}(B), while Fig. \ref{fig:frfs_vf_peaks}(A) shows the estimated peaks via VF overlaid on the average FRF in this frequency band.

The first 30 peaks identified using the VF approach, as shown in Fig. \ref{fig:frfs_vf_peaks}(A), are summarized in Table \ref{tbl3}. In this table, $f_{\mathrm{exp}}$ denotes the modal frequencies obtained experimentally from VF, while $f_{\mathrm{mod}}$ corresponds to the frequencies predicted by the analytical model described in Section~\ref{sec:Section2}. The comparison highlights a good agreement between the experimental and analytical results, especially as the frequency increases. As mentioned earlier, the higher error at lower frequencies is likely due to averaging VEM properties.

\begin{table}[pos=h] 
  \caption{Experimental vs.\ analytical}
  \label{tbl3}

  \begin{tabular*}{\tblwidth}{@{}LLLLLLLL@{}}
  \toprule
  \textbf{Mode} & \textbf{$f_{\text{exp}}$} & \textbf{$f_{\text{mod}}$} & \textbf{Error \%} %
                & \textbf{Mode} & \textbf{$f_{\text{exp}}$} & \textbf{$f_{\text{mod}}$} & \textbf{Error \%}\\
  \midrule
        1 &  12.2193 &  11.5393 & 5.56 & 16 & 765.015  & 813.396  & 6.32\\
        2 &  18.8713 &  18.6817 &  1.00 & 17 & 864.022  & 913.028  & 5.67\\
        3 &  37.6478 &  37.8160 &  0.45 & 18 & 957.644  & 1019.00  & 6.41\\
        4 &  62.2832 &  68.1554 &  9.43 & 19 & 1060.69  & 1132.35  & 6.76\\
        5 &  81.6156 &  90.4775 &  10.9 & 20 & 1187.65  & 1239.70  & 4.38\\
        6 &  122.960 &  125.753 & 2.27 & 21 & 1312.15  & 1367.91  & 4.25\\
        7 &  164.406 &  176.788 & 7.53 & 22 & 1422.73  & 1508.65  & 6.04\\
        8 &  199.131 &  221.947 &  11.5 & 23 & 1554.47  & 1642.08  & 5.64\\
        9 &  254.010 &  268.536 &  5.72 & 24 & 1679.28  & 1779.59  & 5.97\\
       10 &  308.705 &  333.602 &  8.07 & 25 & 1839.98  & 1922.24  & 4.47\\
       11 &  370.582 &  396.780 &  7.07 & 26 & 2015.45  & 2069.59  & 2.69\\
       12 &  441.545 &  467.104 &  5.79 & 27 & 2153.24  & 2241.48  & 4.10\\
       13 &  504.586 &  549.981 &  9.00 & 28 & 2266.92  & 2419.96  & 6.75\\
       14 &  585.686 &  625.127 &  6.73 & 29 & 2440.61  & 2579.44  & 5.69\\
       15 &  676.707 &  711.041 &  5.07 & 30 & 2642.12  & 2744.81  & 3.89\\
  \bottomrule
\end{tabular*}
\end{table}

\section{Results and Discussion}

In order to assess the wave quality, a cost function CF discussed in Ref. \cite{Motaharibidgoli_2023} is used as follows,

\begin{equation}
  CF=\frac{v_{max}-v_{min}}{v_{max}+v_{min}}.
  \label{eq:18}
\end{equation}

When TW is dominant in the response, CF yields lower values, with 0 suggesting a pure TW response. In contrast, higher values indicate a standing wave (SW)-dominant response, with 1 representing a pure SW response. In the following sections, the results obtained using the setup shown in Fig. \ref{fig:EBABHExp} are compared against those yielded by the developed model for validation purposes. Next, a parametric study is conducted using the validated model to examine various parameters affecting the system's response.

\subsection{Model Validation}
\label{sub:MV}

\begin{figure}
    \centering
    \includegraphics[width=1\textwidth]{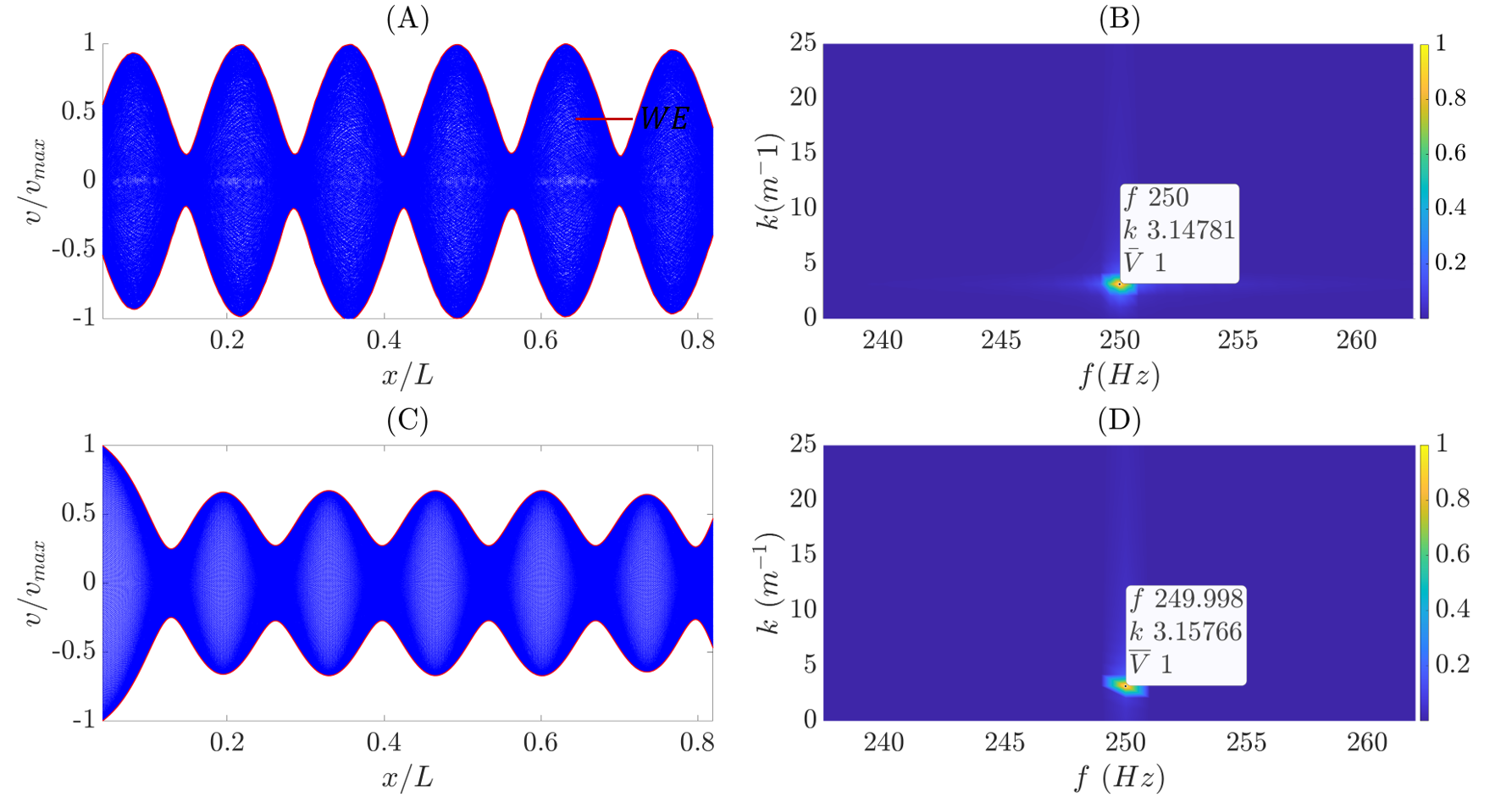}
    \caption{(A) Experimental Wave Envelope, (B) Experimental Response's 2D FFT, (C) Analytical Wave Envelope, and (D) Analytical Response's 2D FFT for the $250~Hz$ case.}
    \label{fig:250hz}
\end{figure}

In this section, the analytical and experimental responses are compared at two excitation frequencies spanning the lower and upper limits of the band of interest. Fig. \ref{fig:250hz} shows the fully developed analytical and experimental time histories across the beam's uniform section ($x=50~mm$ to $x=1000~mm$) as well as their frequency and wave number components obtained by applying a 2D FFT to the responses when the beam is excited at $250~ Hz$. Note that $\overline{V}=v/v_{max}$. As shown in Fig. \ref{fig:250hz}(B) and Fig. \ref{fig:250hz}(D), there is an excellent agreement between the experimental and analytical responses in terms of frequency and wave number components. Despite this, a discrepancy exists in the TW (or SW) contribution in the response, which, as mentioned earlier, is mostly attributed to the use of averaged properties for the damping material, as seen in Fig. \ref{fig:250hz}(A) and Fig. \ref{fig:250hz}(C). The CF values for the experimental and analytical cases are 0.67 and 0.42, respectively. It is expected that the model will become more accurate as the excitation frequency increases.

\begin{figure}
    \centering
    \includegraphics[width=1\textwidth]{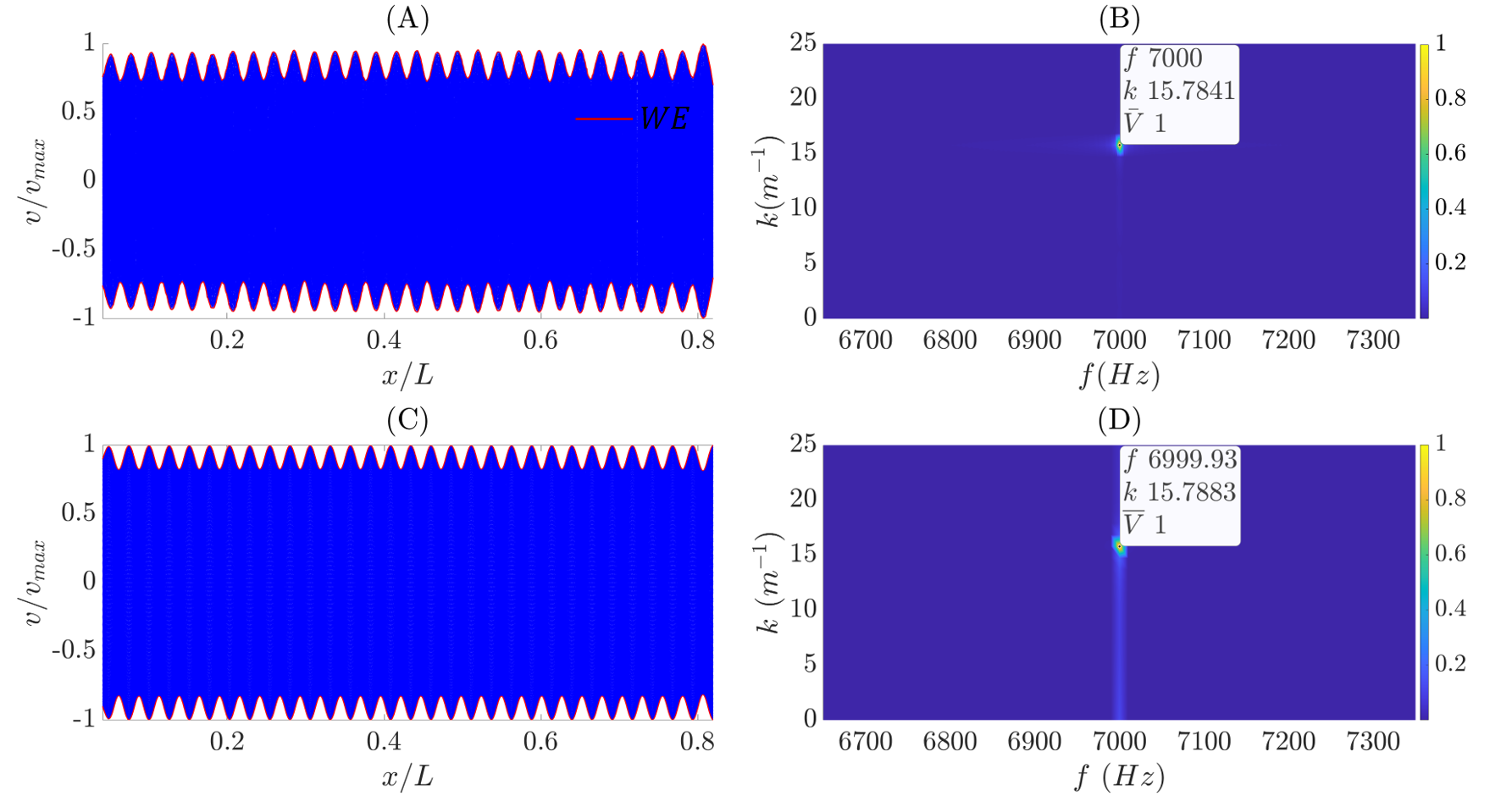}
    \caption{(A) Experimental Wave Envelope, (B) Experimental Response's 2D FFT, (C) Analytical Wave Envelope, and (D) Analytical Response's 2D FFT.}
    \label{fig:7000hz}
\end{figure}

\begin{figure}
    \centering
    \includegraphics[width=0.85\textwidth]{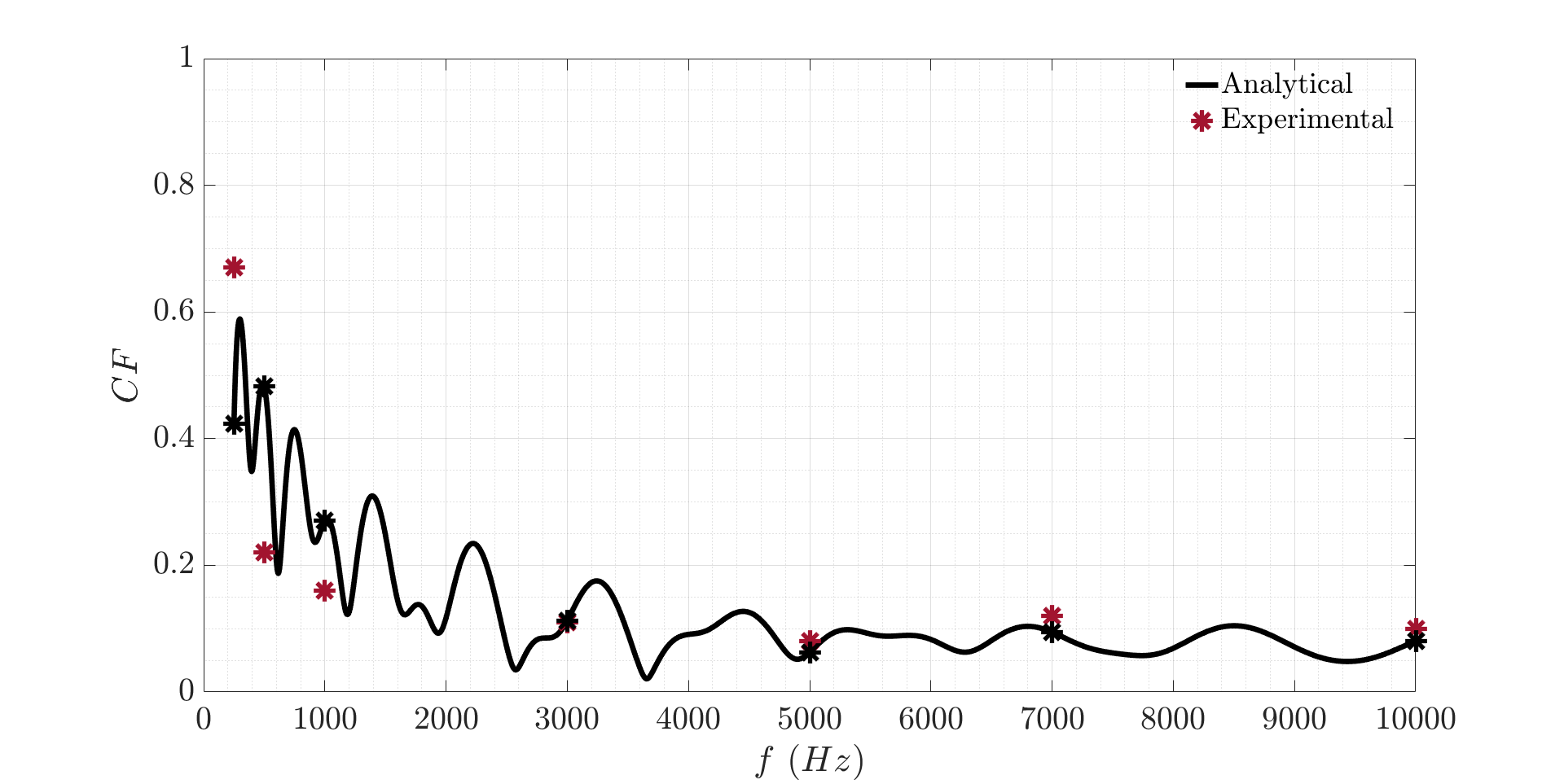}
    \caption{CF vs. Frequency for experimental and analytical cases.}
    \label{fig:cf-f}
\end{figure}

Thus, as shown in the following example, the case where the system is excited at $7000~Hz$ is inspected in Fig. \ref{fig:7000hz}. This reveals excellent agreement between the responses in both frequency and time domains. This is due, as mentioned before and shown in Table \ref{tbl3}, to the assumption that the mean loss and storage moduli yield more accurate results as the frequency increases. The CF values of 0.12 and 0.1 are obtained from the experimental results and the model, respectively. As such, it is concluded that the model closely follows the experiment. The acceptable frequency range for the model’s representation of the experiment can be established from Fig. \ref{fig:cf-f}.

The CF values associated with six different frequencies are calculated following experimental analyses and are depicted in Fig. \ref{fig:cf-f}. Comparing the developed model's response to the aforementioned experimental cases, it can be deduced that the model is increasingly more accurate in the range of $[1~kHz-10~kHz]$. As such, this range will be considered for the parametric study.  

Also, Fig. \ref{fig:cf-f} illustrates the sensitivity of 
$CF$ to excitation frequency. Despite the oscillatory nature of the response, $CF$ exhibits an overall decreasing trend with increasing frequency, consistently observed in both analytical and experimental results. Moreover, the oscillations diminish in amplitude at higher frequencies, further confirming the enhanced effectiveness of ABH treatment at elevated excitation frequencies \cite{austin2022realisation,yang2025study}. This suggests that, despite the general increase in TW's dominance at higher frequencies, near-perfect TWs occur at certain instances. However, for a response to be acceptable in terms of TW quality, it is not required that the CF be zero; higher values may be acceptable depending on the specific application in mind. As such, based on the authors' observation in this specific case study, in frequencies above $2~kHz$ the $CF$ is almost always below $0.2$, suggesting at least about $80\%$ TW contribution in the response, which is significant considering that perfect TW ($CF=0$) is practically impossible in dispersive systems like beams. Thus, ABHs could be considered an effective solution for the broadband TW problem in various applications, including the development of biomechanical models that mimic the vibratory behavior of the BM.

\subsection{Parametric Study}
In the following sections, the effects of three different system parameters, i.e., viscoelastic material loss factor, power-law order, and tapered section length, on the system's response are examined.

\subsubsection{Viscoelastic material loss factor}

\begin{figure}
    \centering
    \includegraphics[width=0.85\textwidth]{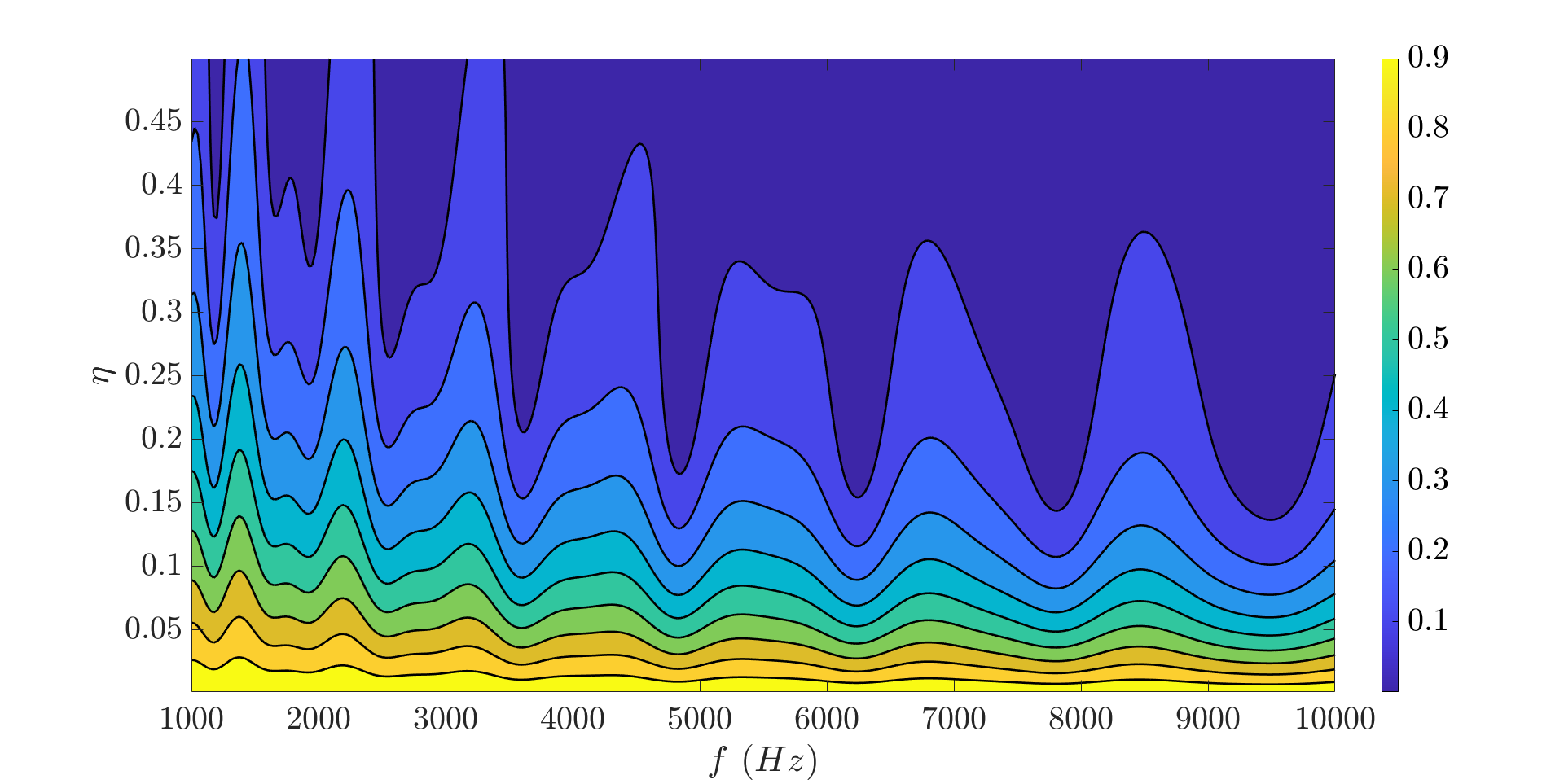}
    \caption{The variation of $CF$ with the VEM loss factor and excitation frequency.}
    \label{fig:LFFCF}
\end{figure}

In this section, the effect of the VEM loss factor on the system's response is investigated. It is assumed that the loss factor varies gradually over the range $0.001$ to $0.5$, spanning the practical range from lightly damped structural materials to highly damped VEMs. Fig. \ref{fig:LFFCF} shows the $ CF$'s trend associated with the VEM's loss factor and excitation frequency variation. 

The result reveals that for a given excitation frequency, increasing the loss factor increases the TW dominance in the response (decreasing the $CF$). The results imply that at higher frequencies, the $CF$ decreases at a higher rate as $\eta$ increases, suggesting that even a smaller amount of added damping could promote TW production in such cases. Furthermore, the oscillatory behavior with frequency variation, which was previously discussed, is also observable here. However, it should be noted that for very low values of $\eta$, no matter how high the frequencies, there wouldn't be a high-quality TW response. The loss factor should be at least $0.1$ to enable a roughly $70\%$ TW response at high frequencies. As the final observation in this section, it can be seen that at the higher end of the $\eta$ range, high-quality TWs can be obtained in frequencies as low as $1000~Hz$.

\subsubsection{Power-law order}

\begin{figure}
    \centering
    \includegraphics[width=0.85\textwidth]{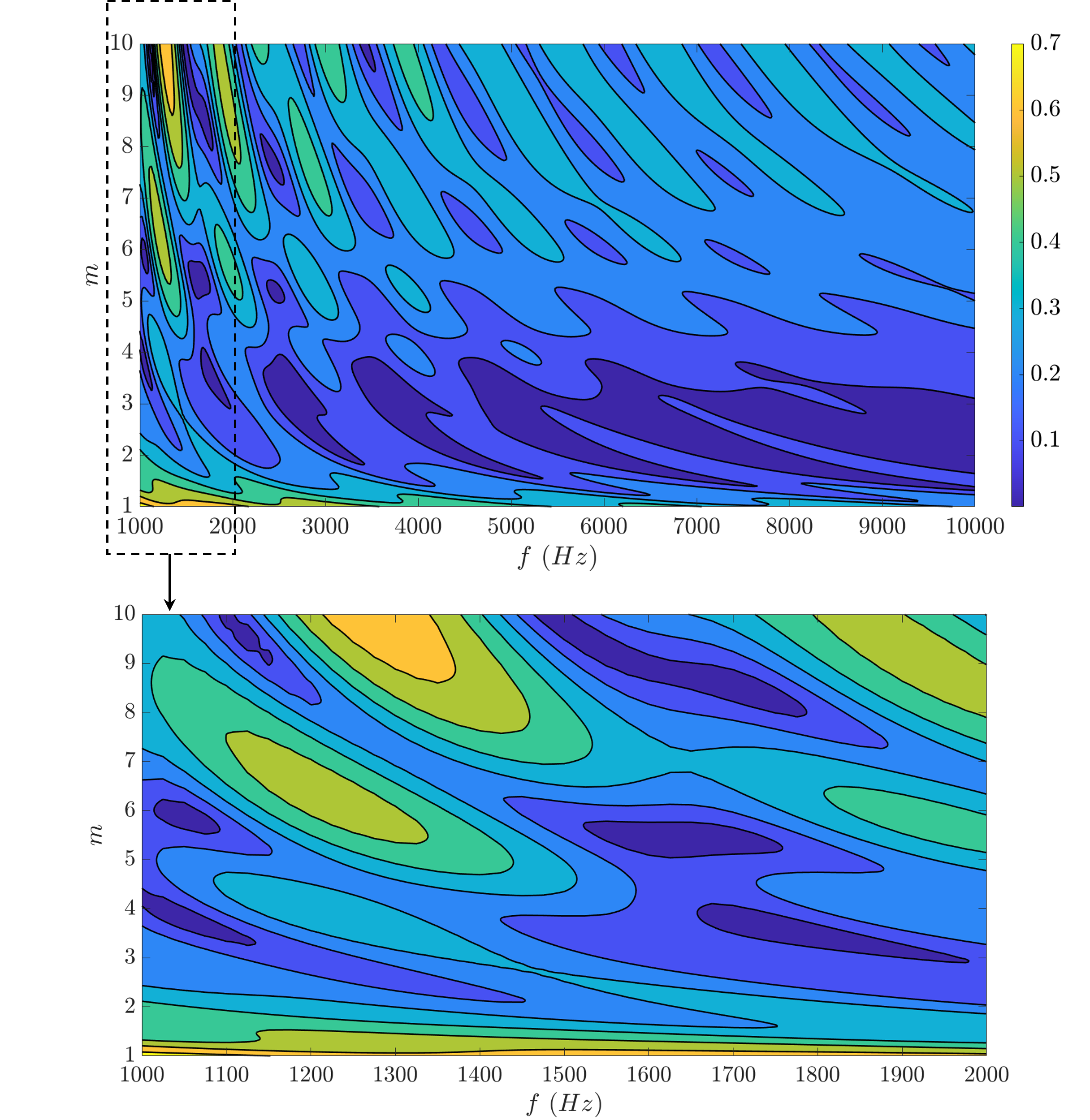}
    \caption{The variation of $CF$ with the power-law order and excitation frequency.}
    \label{fig:mFCF}
\end{figure}

The effects of the power-law order, $m$, in the tapered section are studied and presented in this section. It is assumed that $m$ smoothly transitions from a linear tapering profile, i.e., $m=1$, to a highly nonlinear one, $m=10$. The trend is depicted in Fig. \ref{fig:mFCF}.

Considering a given exponent, it is observed that by increasing the frequency, the $CF$ generally decreases in a similar oscillatory fashion that was observed in Fig. \ref{fig:cf-f}. However, it is worth noting, especially towards the lower excitation frequency end, that $CF$ is generally lower when $2\leq m\leq4$. This confirms the observations reported in the literature \cite{hook2019}. It can also be observed that as $m$ increases, the $CF$ decreases before increasing again, in an oscillatory manner. At frequencies above $4~kHz$ this fluctuation is still within an acceptable $CF$ range, where there is at least about $70\%$ TW contribution response at any given $m$. Finally, it can be deduced from what is presented and discussed in this section that, contingent on choosing an exponent around $3$, high-quality TWs could be obtained in the range of $[1~kHz -10~kHz$] and potentially higher. 

\subsubsection{Tapered section length}

\begin{figure}
    \centering
    \includegraphics[width=0.85\textwidth]{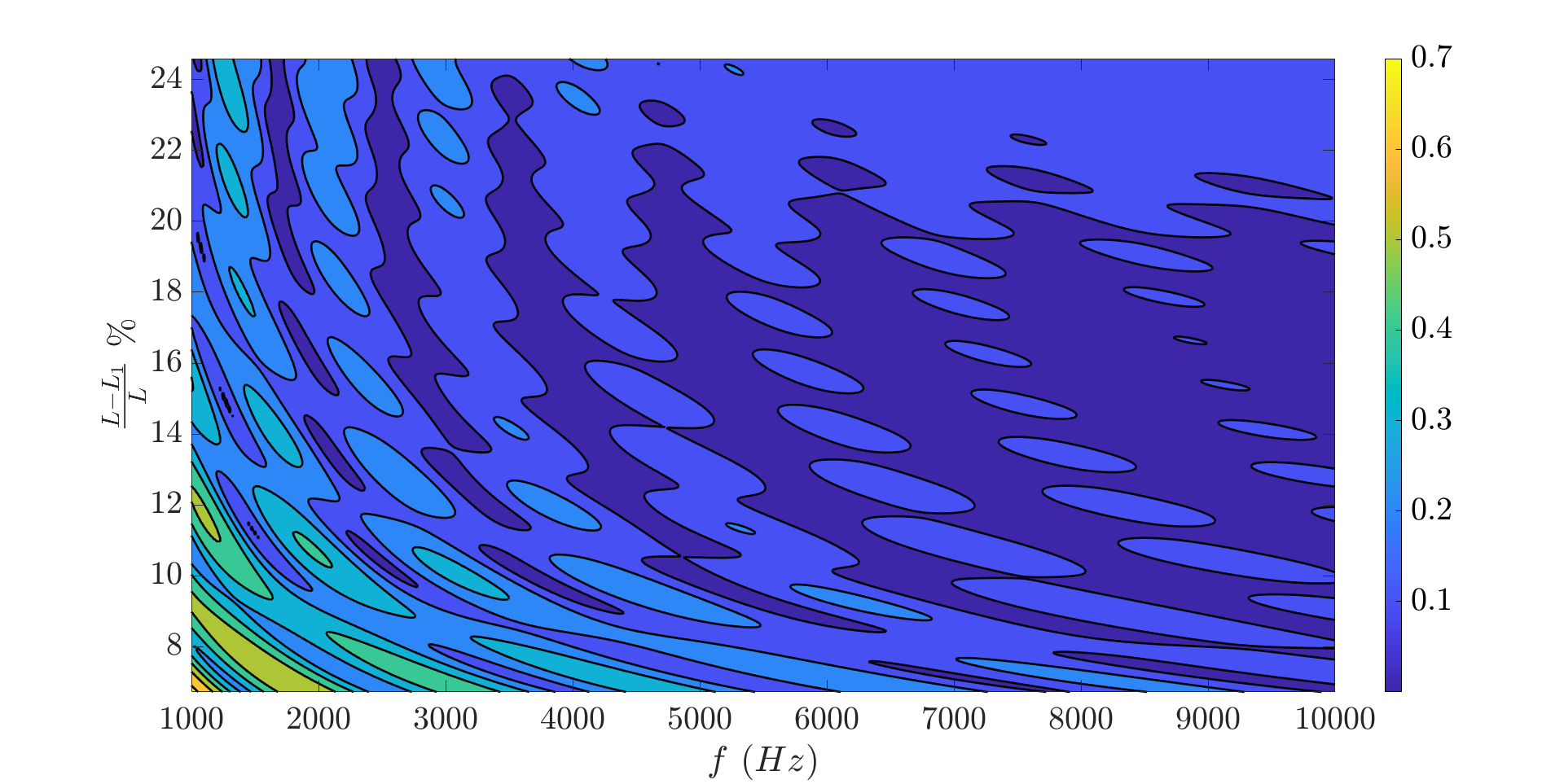}
    \caption{The variation of $CF$ with the ABH section length and excitation frequency.}
    \label{fig:L1FCF}
\end{figure}

This section aims to investigate the effects of varying the ABH termination length in comparison to the total system length. It starts with the case where the tapered section is about $7\%$ of the total length, fully covered by the VEM tape, and gradually increases to roughly $25\%$. The VEM tape length is kept fixed according to Table \ref{tab:matpropEBABH}. The results are presented in Fig. \ref{fig:L1FCF}. It can be observed that at any given excitation frequency, by increasing the tapered length to the overall length ratio, the $CF$ non-monotonically decreases and then increases, with higher rates of change at higher frequencies. This interestingly suggests that when ABH section is roughly $15\%-20\%$ of the total length (Here, this could also be translated to the VEM coverage of around $34\%-45\%$ of the total ABH length.), high-quality TWs could be obtained for the entire band of interest. A pattern similar to what was observed in Fig. \ref{fig:cf-f} is seen here for given ABH lengths. Finally, the higher the ratio, the faster the TW dominates the response as frequencies increase.

\section{Conclusions}

The dynamic response of a beam with an ABH section at its end was studied in this work. A model was developed and experimentally validated over the range $[1~kHz, 10~kHz]$. This model was then used to investigate the effects of the system's parameters on the response. The key findings of this study are listed as follows.

\begin{itemize}
  \item Based on observations reported herein, provided the VEM's loss factor is sufficient, responses with at least about $80\%$ TW contributions within the range of $[1~ kHz-10~kHz]$ are achievable.
  \item The presented results suggest that a second to fourth-order power-law tapering generally leads to better performance, with 3rd order being optimal. Values outside this range would weaken the response in the lower frequencies. However, higher power-law exponents would not significantly affect the response at high frequencies as compared to lower frequencies.
  \item It can be concluded that the VEM layer needs to partially cover the ABH section to improve the response in terms of its TW component for the entire frequency range. Similar to the power-law order, it appears that there is a sweet spot for the VEM coverage, which is around $34\% - 45\%$ of the total length of the ABH. This happens to coincide with the case that the ABH is around $15\% - 20\%$ of the total length, $L$.
  
\end{itemize}

Considering these points, the parameters chosen to develop the experimental setup were mostly optimal as high-quality TWs ($CF\leq0.2$) were sustained over the range $[500~Hz-10~kHz$]. Ultimately, these findings suggest that ABHs are a favorable solution for generating TWs over a wide range of excitation frequencies, should the system be designed considering the findings reported above. 

\section*{Acknowledgements}

AO and SS would like to recognize the support provided by Texas A\&M's J. Mike Walker '66 Department of Mechanical Engineering through the Byron Anderson '54 Fellowship and Graduate Summer Research Grant. SNHS gratefully acknowledges the support of the 2024 Fulbright Visiting Scholar Program. PAT would like to recognize the support provided by the James J. Cain '51 fellowship at Texas A\&M. The authors would like to acknowledge Mr. Trevor Turner and Texas A\&M High Performance Research Computing (HPRC) for providing the advanced computing resources needed for this study. 

\printcredits

\bibliographystyle{model1-num-names}

\bibliography{cas-refs}

\begin{thebibliography}{79}
\expandafter\ifx\csname natexlab\endcsname\relax\def\natexlab#1{#1}\fi
\providecommand{\url}[1]{\texttt{#1}}
\providecommand{\href}[2]{#2}
\providecommand{\path}[1]{#1}
\providecommand{\DOIprefix}{doi:}
\providecommand{\ArXivprefix}{arXiv:}
\providecommand{\URLprefix}{URL: }
\providecommand{\Pubmedprefix}{pmid:}
\providecommand{\doi}[1]{\href{http://dx.doi.org/#1}{\path{#1}}}
\providecommand{\Pubmed}[1]{\href{pmid:#1}{\path{#1}}}
\providecommand{\bibinfo}[2]{#2}
\ifx\xfnm\relax \def\xfnm[#1]{\unskip,\space#1}\fi
\bibitem[{Taylor(1952)}]{taylor1952action}
\bibinfo{author}{G.~I. Taylor},
\newblock \bibinfo{title}{The action of waving cylindrical tails in propelling microscopic organisms},
\newblock \bibinfo{journal}{Proceedings of the Royal Society of London. Series A. Mathematical and Physical Sciences} \bibinfo{volume}{211} (\bibinfo{year}{1952}) \bibinfo{pages}{225--239}.
\bibitem[{Stone and Samuel(1996)}]{stone1996propulsion}
\bibinfo{author}{H.~A. Stone}, \bibinfo{author}{A.~D. Samuel},
\newblock \bibinfo{title}{Propulsion of microorganisms by surface distortions},
\newblock \bibinfo{journal}{Physical review letters} \bibinfo{volume}{77} (\bibinfo{year}{1996}) \bibinfo{pages}{4102}.
\bibitem[{Vig and Wolgemuth(2012)}]{vig2012swimming}
\bibinfo{author}{D.~K. Vig}, \bibinfo{author}{C.~W. Wolgemuth},
\newblock \bibinfo{title}{Swimming dynamics of the lyme disease spirochete},
\newblock \bibinfo{journal}{Physical review letters} \bibinfo{volume}{109} (\bibinfo{year}{2012}) \bibinfo{pages}{218104}.
\bibitem[{Fish(1982)}]{fish1982function}
\bibinfo{author}{F.~E. Fish},
\newblock \bibinfo{title}{Function of the compressed tail of surface swimming muskrats (ondatra zibethicus)},
\newblock \bibinfo{journal}{Journal of Mammalogy} \bibinfo{volume}{63} (\bibinfo{year}{1982}) \bibinfo{pages}{591--597}.
\bibitem[{Long~Jr and Nipper(1996)}]{long1996importance}
\bibinfo{author}{J.~H. Long~Jr}, \bibinfo{author}{K.~S. Nipper},
\newblock \bibinfo{title}{The importance of body stiffness in undulatory propulsion},
\newblock \bibinfo{journal}{American Zoologist} \bibinfo{volume}{36} (\bibinfo{year}{1996}) \bibinfo{pages}{678--694}.
\bibitem[{Marvi et~al.(2013)Marvi, Bridges, and Hu}]{marvi2013snakes}
\bibinfo{author}{H.~Marvi}, \bibinfo{author}{J.~Bridges}, \bibinfo{author}{D.~L. Hu},
\newblock \bibinfo{title}{Snakes mimic earthworms: propulsion using rectilinear travelling waves},
\newblock \bibinfo{journal}{Journal of the Royal Society Interface} \bibinfo{volume}{10} (\bibinfo{year}{2013}) \bibinfo{pages}{20130188}.
\bibitem[{Chong et~al.(2022)Chong, Wang, Erickson, Bergmann, and Goldman}]{chong2022coordinating}
\bibinfo{author}{B.~Chong}, \bibinfo{author}{T.~Wang}, \bibinfo{author}{E.~Erickson}, \bibinfo{author}{P.~J. Bergmann}, \bibinfo{author}{D.~I. Goldman},
\newblock \bibinfo{title}{Coordinating tiny limbs and long bodies: Geometric mechanics of lizard terrestrial swimming},
\newblock \bibinfo{journal}{Proceedings of the National Academy of Sciences} \bibinfo{volume}{119} (\bibinfo{year}{2022}) \bibinfo{pages}{e2118456119}.
\bibitem[{Von~B{\'e}k{\'e}sy(1928)}]{von1928theorie}
\bibinfo{author}{G.~Von~B{\'e}k{\'e}sy}, \bibinfo{title}{Zur Theorie des H{\"o}rens: Die Schwingungsform der Basilarmembran}, \bibinfo{publisher}{{\'E}diteur inconnu}, \bibinfo{year}{1928}.
\bibitem[{Von~B{\'e}k{\'e}sy(1960)}]{von1960experiments}
\bibinfo{author}{G.~Von~B{\'e}k{\'e}sy}, \bibinfo{title}{Experiments in hearing.}, \bibinfo{publisher}{McGraw Hill}, \bibinfo{year}{1960}.
\bibitem[{voN B{\'e}k{\'e}sY(1970)}]{von1970travelling}
\bibinfo{author}{G.~voN B{\'e}k{\'e}sY},
\newblock \bibinfo{title}{Travelling waves as frequency analysers in the cochlea},
\newblock \bibinfo{journal}{Nature} \bibinfo{volume}{225} (\bibinfo{year}{1970}) \bibinfo{pages}{1207--1209}.
\bibitem[{Zwislocki(1980)}]{zwislocki1980theory}
\bibinfo{author}{J.~Zwislocki},
\newblock \bibinfo{title}{Theory of cochlear mechanics},
\newblock \bibinfo{journal}{Hearing Research} \bibinfo{volume}{2} (\bibinfo{year}{1980}) \bibinfo{pages}{171--182}.
\bibitem[{Loh and Ro(2000)}]{loh2000object}
\bibinfo{author}{B.-G. Loh}, \bibinfo{author}{P.~I. Ro},
\newblock \bibinfo{title}{An object transport system using flexural ultrasonic progressive waves generated by two-mode excitation},
\newblock \bibinfo{journal}{IEEE transactions on ultrasonics, ferroelectrics, and frequency control} \bibinfo{volume}{47} (\bibinfo{year}{2000}) \bibinfo{pages}{994--999}.
\bibitem[{Zhao et~al.(2022)Zhao, Mu, Dong, Sun, and Grattan}]{zhao2022requirements}
\bibinfo{author}{J.~Zhao}, \bibinfo{author}{G.~Mu}, \bibinfo{author}{H.~Dong}, \bibinfo{author}{T.~Sun}, \bibinfo{author}{K.~T. Grattan},
\newblock \bibinfo{title}{Requirements for a transportation system based on ultrasonic traveling waves using the measurement of spatial phase difference},
\newblock \bibinfo{journal}{Mechanical Systems and Signal Processing} \bibinfo{volume}{168} (\bibinfo{year}{2022}) \bibinfo{pages}{108708}.
\bibitem[{Rogers et~al.(2023)Rogers, Albakri, and Tarazaga}]{rogers2023directed}
\bibinfo{author}{W.~C. Rogers}, \bibinfo{author}{M.~I. Albakri}, \bibinfo{author}{P.~Tarazaga},
\newblock \bibinfo{title}{Directed particle motion driven by superimposed two-dimensional traveling waves},
\newblock in: \bibinfo{booktitle}{Smart Materials, Adaptive Structures and Intelligent Systems}, volume \bibinfo{volume}{87523}, \bibinfo{organization}{American Society of Mechanical Engineers}, \bibinfo{year}{2023}, p. \bibinfo{pages}{V001T06A003}.
\bibitem[{Mansouri et~al.(2024)Mansouri, Hsiao-Wecksler, and Krishnan}]{mansouri2024toward}
\bibinfo{author}{M.~Mansouri}, \bibinfo{author}{E.~T. Hsiao-Wecksler}, \bibinfo{author}{G.~Krishnan},
\newblock \bibinfo{title}{Toward design guidelines for multidirectional patient transfer on a bed surface using traveling waves},
\newblock \bibinfo{journal}{Journal of Mechanisms and Robotics} \bibinfo{volume}{16} (\bibinfo{year}{2024}) \bibinfo{pages}{074501}.
\bibitem[{Zhu et~al.(2012)Zhu, Gao, and He}]{zhu2012piezoelectric}
\bibinfo{author}{J.~Zhu}, \bibinfo{author}{C.~Gao}, \bibinfo{author}{L.~He},
\newblock \bibinfo{title}{Piezoelectric-based crack detection techniques of concrete structures: Experimental study},
\newblock \bibinfo{journal}{Journal of Wuhan University of Technology-Mater. Sci. Ed.} \bibinfo{volume}{27} (\bibinfo{year}{2012}) \bibinfo{pages}{346--352}.
\bibitem[{Lugovtsova et~al.(2019)Lugovtsova, Bulling, Boller, and Prager}]{lugovtsova2019analysis}
\bibinfo{author}{Y.~Lugovtsova}, \bibinfo{author}{J.~Bulling}, \bibinfo{author}{C.~Boller}, \bibinfo{author}{J.~Prager},
\newblock \bibinfo{title}{Analysis of guided wave propagation in a multi-layered structure in view of structural health monitoring},
\newblock \bibinfo{journal}{Applied Sciences} \bibinfo{volume}{9} (\bibinfo{year}{2019}) \bibinfo{pages}{4600}.
\bibitem[{Musgrave and Tarazaga(2019)}]{musgrave2019turbulent}
\bibinfo{author}{P.~F. Musgrave}, \bibinfo{author}{P.~A. Tarazaga},
\newblock \bibinfo{title}{Turbulent boundary layer over a piezoelectrically excited traveling wave surface},
\newblock in: \bibinfo{booktitle}{AIAA Scitech 2019 Forum}, \bibinfo{year}{2019}, p. \bibinfo{pages}{1354}.
\bibitem[{Olivett et~al.(2021)Olivett, Corrao, and Karami}]{olivett2021flow}
\bibinfo{author}{A.~Olivett}, \bibinfo{author}{P.~Corrao}, \bibinfo{author}{M.~A. Karami},
\newblock \bibinfo{title}{Flow control and separation delay in morphing wing aircraft using traveling wave actuation},
\newblock \bibinfo{journal}{Smart Materials and Structures} \bibinfo{volume}{30} (\bibinfo{year}{2021}) \bibinfo{pages}{025028}.
\bibitem[{Zhang et~al.(2024)Zhang, Lv, Zhao, Fan, Xie, Shan, and Li}]{zhang2024achieving}
\bibinfo{author}{L.~Zhang}, \bibinfo{author}{M.~Lv}, \bibinfo{author}{X.~Zhao}, \bibinfo{author}{H.~Fan}, \bibinfo{author}{T.~Xie}, \bibinfo{author}{X.~Shan}, \bibinfo{author}{K.~Li},
\newblock \bibinfo{title}{Achieving travelling wave drag reduction by micro piezoelectric actuator},
\newblock \bibinfo{journal}{International Journal of Mechanical Sciences} \bibinfo{volume}{275} (\bibinfo{year}{2024}) \bibinfo{pages}{109326}.
\bibitem[{Ogunka and Borazjani(2025)}]{ogunka2025simulations}
\bibinfo{author}{U.~E. Ogunka}, \bibinfo{author}{I.~Borazjani},
\newblock \bibinfo{title}{Simulations of wave generation over a flexible structure with piezoelectric actuators for flow control},
\newblock in: \bibinfo{booktitle}{AIAA SCITECH 2025 Forum}, \bibinfo{year}{2025}, p. \bibinfo{pages}{1688}.
\bibitem[{Chen et~al.(2007)Chen, Ma, Wang, Li, and Duan}]{chen2007design}
\bibinfo{author}{L.~Chen}, \bibinfo{author}{S.~Ma}, \bibinfo{author}{Y.~Wang}, \bibinfo{author}{B.~Li}, \bibinfo{author}{D.~Duan},
\newblock \bibinfo{title}{Design and modelling of a snake robot in traveling wave locomotion},
\newblock \bibinfo{journal}{Mechanism and Machine Theory} \bibinfo{volume}{42} (\bibinfo{year}{2007}) \bibinfo{pages}{1632--1642}.
\bibitem[{Qi et~al.(2020)Qi, Shi, Pinto, and Tan}]{qi2020novel}
\bibinfo{author}{X.~Qi}, \bibinfo{author}{H.~Shi}, \bibinfo{author}{T.~Pinto}, \bibinfo{author}{X.~Tan},
\newblock \bibinfo{title}{A novel pneumatic soft snake robot using traveling-wave locomotion in constrained environments},
\newblock \bibinfo{journal}{IEEE Robotics and Automation Letters} \bibinfo{volume}{5} (\bibinfo{year}{2020}) \bibinfo{pages}{1610--1617}.
\bibitem[{Zhu et~al.(2024)Zhu, Li, Wu, and Li}]{zhu2024double}
\bibinfo{author}{B.~Zhu}, \bibinfo{author}{C.~Li}, \bibinfo{author}{Z.~Wu}, \bibinfo{author}{Y.~Li},
\newblock \bibinfo{title}{A double-beam piezoelectric robot based on the principle of two-mode excitation},
\newblock \bibinfo{journal}{Sensors and Actuators A: Physical} \bibinfo{volume}{369} (\bibinfo{year}{2024}) \bibinfo{pages}{115154}.
\bibitem[{Ji and Song(2024)}]{ji2024bionic}
\bibinfo{author}{Q.~Ji}, \bibinfo{author}{A.~Song},
\newblock \bibinfo{title}{Bionic snail robot enhanced by poroelastic foams crawls using direct and retrograde waves},
\newblock \bibinfo{journal}{Soft Robotics} \bibinfo{volume}{11} (\bibinfo{year}{2024}) \bibinfo{pages}{453--463}.
\bibitem[{Malladi et~al.(2017)Malladi, Albakri, Musgrave, and Tarazaga}]{malladi2017investigation}
\bibinfo{author}{V.~V.~S. Malladi}, \bibinfo{author}{M.~Albakri}, \bibinfo{author}{P.~Musgrave}, \bibinfo{author}{P.~A. Tarazaga},
\newblock \bibinfo{title}{Investigation of propulsive characteristics due to traveling waves in continuous finite media},
\newblock in: \bibinfo{booktitle}{Bioinspiration, Biomimetics, and Bioreplication 2017}, volume \bibinfo{volume}{10162}, \bibinfo{organization}{SPIE}, \bibinfo{year}{2017}, pp. \bibinfo{pages}{112--119}.
\bibitem[{Musgrave(2021)}]{musgrave2021electro}
\bibinfo{author}{P.~F. Musgrave},
\newblock \bibinfo{title}{Electro-hydro-elastic modeling of structure-borne traveling waves and their application to aquatic swimming motions},
\newblock \bibinfo{journal}{Journal of Fluids and Structures} \bibinfo{volume}{102} (\bibinfo{year}{2021}) \bibinfo{pages}{103230}.
\bibitem[{Gupta and Malladi(2024)}]{gupta2024utilizing}
\bibinfo{author}{S.~Gupta}, \bibinfo{author}{V.~V.~S. Malladi},
\newblock \bibinfo{title}{Utilizing steady-state traveling waves in a quiescent water environment for particle propulsion},
\newblock in: \bibinfo{booktitle}{Smart Materials, Adaptive Structures and Intelligent Systems}, volume \bibinfo{volume}{88322}, \bibinfo{organization}{American Society of Mechanical Engineers}, \bibinfo{year}{2024}, p. \bibinfo{pages}{V001T03A008}.
\bibitem[{Syuhri et~al.(2024)Syuhri, Zare-Behtash, and Cammarano}]{syuhri2024travelling}
\bibinfo{author}{S.~N. Syuhri}, \bibinfo{author}{H.~Zare-Behtash}, \bibinfo{author}{A.~Cammarano},
\newblock \bibinfo{title}{Travelling waves in beam-like structures submerged in water},
\newblock \bibinfo{journal}{International Journal of Mechanical Sciences} \bibinfo{volume}{283} (\bibinfo{year}{2024}) \bibinfo{pages}{109623}.
\bibitem[{Hess and Musgrave(2024)}]{hess2024continuum}
\bibinfo{author}{I.~Hess}, \bibinfo{author}{P.~Musgrave},
\newblock \bibinfo{title}{A continuum soft robotic trout with embedded hasel actuators: design, fabrication, and swimming kinematics},
\newblock \bibinfo{journal}{Smart Materials and Structures} \bibinfo{volume}{33} (\bibinfo{year}{2024}) \bibinfo{pages}{105043}.
\bibitem[{Hernandez et~al.(2013)Hernandez, Bernard, and Razek}]{hernandez2013design}
\bibinfo{author}{C.~Hernandez}, \bibinfo{author}{Y.~Bernard}, \bibinfo{author}{A.~Razek},
\newblock \bibinfo{title}{Design and manufacturing of a piezoelectric traveling-wave pumping device},
\newblock \bibinfo{journal}{IEEE Transactions on ultrasonics, ferroelectrics, and frequency control} \bibinfo{volume}{60} (\bibinfo{year}{2013}) \bibinfo{pages}{1949--1956}.
\bibitem[{Rajendran et~al.(2023)Rajendran, Gosavi, and Malladi}]{rajendran2023novel}
\bibinfo{author}{K.~Rajendran}, \bibinfo{author}{H.~Gosavi}, \bibinfo{author}{V.~V.~S. Malladi},
\newblock \bibinfo{title}{Novel pumping mechanism for heat sinks with fluid medium using steady state traveling waves},
\newblock in: \bibinfo{booktitle}{Smart Materials, Adaptive Structures and Intelligent Systems}, volume \bibinfo{volume}{87523}, \bibinfo{organization}{American Society of Mechanical Engineers}, \bibinfo{year}{2023}, p. \bibinfo{pages}{V001T06A014}.
\bibitem[{Alajlouni et~al.(2018)Alajlouni, Albakri, and Tarazaga}]{alajlouni2018impact}
\bibinfo{author}{S.~Alajlouni}, \bibinfo{author}{M.~Albakri}, \bibinfo{author}{P.~Tarazaga},
\newblock \bibinfo{title}{Impact localization in dispersive waveguides based on energy-attenuation of waves with the traveled distance},
\newblock \bibinfo{journal}{Mechanical Systems and Signal Processing} \bibinfo{volume}{105} (\bibinfo{year}{2018}) \bibinfo{pages}{361--376}.
\bibitem[{Alajlouni and Tarazaga(2020)}]{alajlouni2020passive}
\bibinfo{author}{S.~Alajlouni}, \bibinfo{author}{P.~Tarazaga},
\newblock \bibinfo{title}{A passive energy-based method for footstep impact localization, using an underfloor accelerometer sensor network with kalman filtering},
\newblock \bibinfo{journal}{Journal of Vibration and Control} \bibinfo{volume}{26} (\bibinfo{year}{2020}) \bibinfo{pages}{941--951}.
\bibitem[{Alajlouni et~al.(2022)Alajlouni, Baker, and Tarazaga}]{alajlouni2022maximum}
\bibinfo{author}{S.~Alajlouni}, \bibinfo{author}{J.~Baker}, \bibinfo{author}{P.~Tarazaga},
\newblock \bibinfo{title}{Maximum likelihood estimation for passive energy-based footstep localization},
\newblock \bibinfo{journal}{Mechanical Systems and Signal Processing} \bibinfo{volume}{163} (\bibinfo{year}{2022}) \bibinfo{pages}{108158}.
\bibitem[{Malladi et~al.(2015)Malladi, Avirovik, Priya, and Tarazaga}]{malladi2015characterization}
\bibinfo{author}{V.~Malladi}, \bibinfo{author}{D.~Avirovik}, \bibinfo{author}{S.~Priya}, \bibinfo{author}{P.~Tarazaga},
\newblock \bibinfo{title}{Characterization and representation of mechanical waves generated in piezo-electric augmented beams},
\newblock \bibinfo{journal}{Smart Materials and Structures} \bibinfo{volume}{24} (\bibinfo{year}{2015}) \bibinfo{pages}{105026}.
\bibitem[{Avirovik et~al.(2016)Avirovik, Malladi, Priya, and Tarazaga}]{avirovik2016theoretical}
\bibinfo{author}{D.~Avirovik}, \bibinfo{author}{V.~S. Malladi}, \bibinfo{author}{S.~Priya}, \bibinfo{author}{P.~A. Tarazaga},
\newblock \bibinfo{title}{Theoretical and experimental correlation of mechanical wave formation on beams},
\newblock \bibinfo{journal}{Journal of Intelligent Material Systems and Structures} \bibinfo{volume}{27} (\bibinfo{year}{2016}) \bibinfo{pages}{1939--1948}.
\bibitem[{Anakok et~al.(2022)Anakok, Davaria, Tarazaga, and Malladi}]{anakok2022study}
\bibinfo{author}{I.~Anakok}, \bibinfo{author}{S.~Davaria}, \bibinfo{author}{P.~A. Tarazaga}, \bibinfo{author}{V.~V.~S. Malladi},
\newblock \bibinfo{title}{A study on steady-state traveling waves in one-dimensional non-dispersive finite media},
\newblock \bibinfo{journal}{Journal of Sound and Vibration} \bibinfo{volume}{528} (\bibinfo{year}{2022}) \bibinfo{pages}{116907}.
\bibitem[{Zhu et~al.(2021)Zhu, Li, Wang, and Zhang}]{zhu2021design}
\bibinfo{author}{B.~Zhu}, \bibinfo{author}{C.~Li}, \bibinfo{author}{R.~Wang}, \bibinfo{author}{C.~Zhang},
\newblock \bibinfo{title}{Design and investigation of a new piezoelectric beam transportation device based on two-mode excitation},
\newblock \bibinfo{journal}{Smart Materials and Structures} \bibinfo{volume}{30} (\bibinfo{year}{2021}) \bibinfo{pages}{115012}.
\bibitem[{Cheng et~al.(2025)Cheng, McFarland, Hua, Lu, and Tan}]{cheng2025traveling}
\bibinfo{author}{X.~Cheng}, \bibinfo{author}{D.~M. McFarland}, \bibinfo{author}{X.~Hua}, \bibinfo{author}{H.~Lu}, \bibinfo{author}{C.~A. Tan},
\newblock \bibinfo{title}{Traveling waves in tensioned euler-bernoulli beams with viscoelastic boundary conditions},
\newblock \bibinfo{journal}{International Journal of Mechanical Sciences} \bibinfo{volume}{294} (\bibinfo{year}{2025}) \bibinfo{pages}{110248}.
\bibitem[{Musgrave et~al.(2021)Musgrave, Albakri, and Phoenix}]{musgrave2021guidelines}
\bibinfo{author}{P.~F. Musgrave}, \bibinfo{author}{M.~I. Albakri}, \bibinfo{author}{A.~A. Phoenix},
\newblock \bibinfo{title}{Guidelines and procedure for tailoring high-performance, steady-state traveling waves for propulsion and solid-state motion},
\newblock \bibinfo{journal}{Smart Materials and Structures} \bibinfo{volume}{30} (\bibinfo{year}{2021}) \bibinfo{pages}{025013}.
\bibitem[{Rogers et~al.(2024)Rogers, Soroor, Turner, Albakri, and Tarazaga}]{rogers2024experimental}
\bibinfo{author}{W.~C. Rogers}, \bibinfo{author}{A.~O. Soroor}, \bibinfo{author}{T.~C. Turner}, \bibinfo{author}{M.~I. Albakri}, \bibinfo{author}{P.~Tarazaga},
\newblock \bibinfo{title}{Experimental demonstration of superimposed orthogonal two-dimensional structure-borne traveling waves},
\newblock in: \bibinfo{booktitle}{IMAC, A Conference and Exposition on Structural Dynamics}, \bibinfo{organization}{Springer}, \bibinfo{year}{2024}, pp. \bibinfo{pages}{99--108}.
\bibitem[{Rogers and Albakri(2024)}]{rogers2024estimation}
\bibinfo{author}{W.~C. Rogers}, \bibinfo{author}{M.~I. Albakri},
\newblock \bibinfo{title}{Estimation, tuning, and evaluation of propagation direction behavior in superimposed orthogonal two-dimensional structure-borne traveling waves},
\newblock \bibinfo{journal}{Smart Materials and Structures} \bibinfo{volume}{33} (\bibinfo{year}{2024}) \bibinfo{pages}{125027}.
\bibitem[{Phoenix et~al.(2015)Phoenix, Malladi, and Tarazaga}]{phoenix2015traveling}
\bibinfo{author}{A.~Phoenix}, \bibinfo{author}{V.~S. Malladi}, \bibinfo{author}{P.~A. Tarazaga},
\newblock \bibinfo{title}{Traveling wave phenomenon through piezoelectric actuation of a free-free cylindrical tube},
\newblock in: \bibinfo{booktitle}{Smart Materials, Adaptive Structures and Intelligent Systems}, volume \bibinfo{volume}{57304}, \bibinfo{organization}{American Society of Mechanical Engineers}, \bibinfo{year}{2015}, p. \bibinfo{pages}{V002T04A018}.
\bibitem[{Motaharibidgoli et~al.(2023)Motaharibidgoli, Davaria, Sriram~Malladi, and Tarazaga}]{Motaharibidgoli_2023}
\bibinfo{author}{S.~Motaharibidgoli}, \bibinfo{author}{S.~Davaria}, \bibinfo{author}{V.~V. Sriram~Malladi}, \bibinfo{author}{P.~A. Tarazaga},
\newblock \bibinfo{title}{Developing coexisting traveling and standing waves in euler-bernoulli beams using a single-point excitation and a spring-damper system},
\newblock \bibinfo{journal}{Journal of Sound and Vibration} \bibinfo{volume}{556} (\bibinfo{year}{2023}) \bibinfo{pages}{117728}.
\bibitem[{Blanchard et~al.(2015)Blanchard, McFarland, Bergman, and Vakakis}]{blanchard2015damping}
\bibinfo{author}{A.~Blanchard}, \bibinfo{author}{D.~McFarland}, \bibinfo{author}{L.~Bergman}, \bibinfo{author}{A.~Vakakis},
\newblock \bibinfo{title}{Damping-induced interplay between vibrations and waves in a forced non-dispersive elastic continuum with asymmetrically placed local attachments},
\newblock \bibinfo{journal}{Proceedings of the Royal Society A: Mathematical, Physical and Engineering Sciences} \bibinfo{volume}{471} (\bibinfo{year}{2015}) \bibinfo{pages}{20140402}.
\bibitem[{Cheng et~al.(2017)Cheng, Blanchard, Tan, Lu, Bergman, McFarland, and Vakakis}]{cheng2017separation}
\bibinfo{author}{X.~Cheng}, \bibinfo{author}{A.~Blanchard}, \bibinfo{author}{C.~A. Tan}, \bibinfo{author}{H.~Lu}, \bibinfo{author}{L.~A. Bergman}, \bibinfo{author}{D.~M. McFarland}, \bibinfo{author}{A.~F. Vakakis},
\newblock \bibinfo{title}{Separation of traveling and standing waves in a finite dispersive string with partial or continuous viscoelastic foundation},
\newblock \bibinfo{journal}{Journal of Sound and Vibration} \bibinfo{volume}{411} (\bibinfo{year}{2017}) \bibinfo{pages}{193--209}.
\bibitem[{Cheng et~al.(2019)Cheng, Bergman, McFarland, Tan, Vakakis, and Lu}]{cheng2019co}
\bibinfo{author}{X.~Cheng}, \bibinfo{author}{L.~A. Bergman}, \bibinfo{author}{D.~M. McFarland}, \bibinfo{author}{C.~A. Tan}, \bibinfo{author}{A.~F. Vakakis}, \bibinfo{author}{H.~Lu},
\newblock \bibinfo{title}{Co-existing complexity-induced traveling wave transmission and vibration localization in euler-bernoulli beams},
\newblock \bibinfo{journal}{Journal of Sound and Vibration} \bibinfo{volume}{458} (\bibinfo{year}{2019}) \bibinfo{pages}{22--43}.
\bibitem[{Gupta and Malladi(2024)}]{gupta2024exploring}
\bibinfo{author}{S.~Gupta}, \bibinfo{author}{V.~V.~S. Malladi},
\newblock \bibinfo{title}{Exploring the relationship between cost function of hybrid traveling waves and structural reflection coefficient adapted from acoustics},
\newblock in: \bibinfo{booktitle}{Smart Materials, Adaptive Structures and Intelligent Systems}, volume \bibinfo{volume}{88322}, \bibinfo{organization}{American Society of Mechanical Engineers}, \bibinfo{year}{2024}, p. \bibinfo{pages}{V001T03A007}.
\bibitem[{Xiao et~al.(2017)Xiao, Blanchard, Zhang, Lu, Michael~McFarland, Vakakis, and Bergman}]{xiao2017separation}
\bibinfo{author}{Y.~Xiao}, \bibinfo{author}{A.~Blanchard}, \bibinfo{author}{Y.~Zhang}, \bibinfo{author}{H.~Lu}, \bibinfo{author}{D.~Michael~McFarland}, \bibinfo{author}{A.~F. Vakakis}, \bibinfo{author}{L.~A. Bergman},
\newblock \bibinfo{title}{Separation of traveling and standing waves in a rigid-walled circular duct containing an intermediate impedance discontinuity},
\newblock \bibinfo{journal}{Journal of Vibration and Acoustics} \bibinfo{volume}{139} (\bibinfo{year}{2017}) \bibinfo{pages}{061001}.
\bibitem[{Tanaka and Kikushima(1991)}]{tanaka1991active}
\bibinfo{author}{N.~Tanaka}, \bibinfo{author}{Y.~Kikushima},
\newblock \bibinfo{title}{Active wave control of a flexible beam: proposition of the active sink method},
\newblock \bibinfo{journal}{JSME international journal. Ser. 3, Vibration, control engineering, engineering for industry} \bibinfo{volume}{34} (\bibinfo{year}{1991}) \bibinfo{pages}{159--167}.
\bibitem[{Gabai and Bucher(2008)}]{gabai2008generating}
\bibinfo{author}{R.~Gabai}, \bibinfo{author}{I.~Bucher},
\newblock \bibinfo{title}{Generating traveling vibration waves in finite structures},
\newblock in: \bibinfo{booktitle}{Engineering Systems Design and Analysis}, volume \bibinfo{volume}{48364}, \bibinfo{year}{2008}, pp. \bibinfo{pages}{761--770}.
\bibitem[{Omidi~Soroor and Tarazaga(2023)}]{omidi2023investigation}
\bibinfo{author}{A.~Omidi~Soroor}, \bibinfo{author}{P.~A. Tarazaga},
\newblock \bibinfo{title}{An investigation on the effectiveness of cross-sectional tapering for broadband non-reflective traveling waves generation in beams with passive discontinuities},
\newblock in: \bibinfo{booktitle}{Smart Materials, Adaptive Structures and Intelligent Systems}, volume \bibinfo{volume}{87523}, \bibinfo{organization}{American Society of Mechanical Engineers}, \bibinfo{year}{2023}, p. \bibinfo{pages}{V001T06A002}.
\bibitem[{Soroor and Tarazaga(2025)}]{soroor2025non}
\bibinfo{author}{A.~O. Soroor}, \bibinfo{author}{P.~A. Tarazaga},
\newblock \bibinfo{title}{Non-reflective traveling waves in finite thin beams: A parametric study},
\newblock \bibinfo{journal}{Thin-Walled Structures} \bibinfo{volume}{208} (\bibinfo{year}{2025}) \bibinfo{pages}{112839}.
\bibitem[{Gautier and Krylov(2020)}]{gautier2020recent}
\bibinfo{author}{F.~Gautier}, \bibinfo{author}{V.~V. Krylov},
\newblock \bibinfo{title}{Recent advances in acoustic black hole research},
\newblock \bibinfo{journal}{J. Sound Vib.} \bibinfo{volume}{476} (\bibinfo{year}{2020}) \bibinfo{pages}{115335}.
\bibitem[{Mironov(1988)}]{mironov1988propagation}
\bibinfo{author}{M.~Mironov},
\newblock \bibinfo{title}{Propagation of a flexural wave in a plate whose thickness decreases smoothly to zero in a finite interval},
\newblock \bibinfo{journal}{Sov. Phys. Acoust} \bibinfo{volume}{34} (\bibinfo{year}{1988}) \bibinfo{pages}{318--319}.
\bibitem[{Krylov(2004)}]{krylov2004new}
\bibinfo{author}{V.~V. Krylov},
\newblock \bibinfo{title}{New type of vibration dampers utilising the effect of acoustic'black holes'},
\newblock \bibinfo{journal}{Acta Acustica united with Acustica} \bibinfo{volume}{90} (\bibinfo{year}{2004}) \bibinfo{pages}{830--837}.
\bibitem[{Krylov and Winward(2007)}]{krylov2007experimental}
\bibinfo{author}{V.~V. Krylov}, \bibinfo{author}{R.~Winward},
\newblock \bibinfo{title}{Experimental investigation of the acoustic black hole effect for flexural waves in tapered plates},
\newblock \bibinfo{journal}{Journal of Sound and Vibration} \bibinfo{volume}{300} (\bibinfo{year}{2007}) \bibinfo{pages}{43--49}.
\bibitem[{Bayod(2011)}]{bayod2011application}
\bibinfo{author}{J.~J. Bayod},
\newblock \bibinfo{title}{Application of elastic wedge for vibration damping of turbine blade},
\newblock \bibinfo{journal}{Journal of System Design and Dynamics} \bibinfo{volume}{5} (\bibinfo{year}{2011}) \bibinfo{pages}{1167--1175}.
\bibitem[{Bowyer and Krylov(2014)}]{bowyer2014damping}
\bibinfo{author}{E.~Bowyer}, \bibinfo{author}{V.~V. Krylov},
\newblock \bibinfo{title}{Damping of flexural vibrations in turbofan blades using the acoustic black hole effect},
\newblock \bibinfo{journal}{Applied Acoustics} \bibinfo{volume}{76} (\bibinfo{year}{2014}) \bibinfo{pages}{359--365}.
\bibitem[{Zhang et~al.(2024)Zhang, Ding, Wu, Ma, and Deng}]{zhang2024vibroacoustic}
\bibinfo{author}{S.~Zhang}, \bibinfo{author}{L.~Ding}, \bibinfo{author}{X.~Wu}, \bibinfo{author}{Y.~Ma}, \bibinfo{author}{Z.~Deng},
\newblock \bibinfo{title}{Vibroacoustic suppression of sandwich plates with imperfect acoustic black hole},
\newblock \bibinfo{journal}{International Journal of Mechanical Sciences} \bibinfo{volume}{283} (\bibinfo{year}{2024}) \bibinfo{pages}{109690}.
\bibitem[{Zhao et~al.(2014)Zhao, Conlon, and Semperlotti}]{zhao2014broadband}
\bibinfo{author}{L.~Zhao}, \bibinfo{author}{S.~C. Conlon}, \bibinfo{author}{F.~Semperlotti},
\newblock \bibinfo{title}{Broadband energy harvesting using acoustic black hole structural tailoring},
\newblock \bibinfo{journal}{Smart materials and structures} \bibinfo{volume}{23} (\bibinfo{year}{2014}) \bibinfo{pages}{065021}.
\bibitem[{Ji et~al.(2019)Ji, Liang, Qiu, Cheng, and Wu}]{ji2019enhancement}
\bibinfo{author}{H.~Ji}, \bibinfo{author}{Y.~Liang}, \bibinfo{author}{J.~Qiu}, \bibinfo{author}{L.~Cheng}, \bibinfo{author}{Y.~Wu},
\newblock \bibinfo{title}{Enhancement of vibration based energy harvesting using compound acoustic black holes},
\newblock \bibinfo{journal}{Mechanical Systems and Signal Processing} \bibinfo{volume}{132} (\bibinfo{year}{2019}) \bibinfo{pages}{441--456}.
\bibitem[{Du et~al.(2024)Du, Xiang, and Qiu}]{du2024semi}
\bibinfo{author}{W.~Du}, \bibinfo{author}{Z.~Xiang}, \bibinfo{author}{X.~Qiu},
\newblock \bibinfo{title}{A semi-analytical electromechanical model for energy harvesting of plate with acoustic black hole indentations},
\newblock \bibinfo{journal}{Thin-Walled Structures} \bibinfo{volume}{203} (\bibinfo{year}{2024}) \bibinfo{pages}{112235}.
\bibitem[{Bezan{\c{c}}on et~al.(2024)Bezan{\c{c}}on, Doutres, Umnova, Leclaire, and Dupont}]{bezanccon2024thin}
\bibinfo{author}{G.~Bezan{\c{c}}on}, \bibinfo{author}{O.~Doutres}, \bibinfo{author}{O.~Umnova}, \bibinfo{author}{P.~Leclaire}, \bibinfo{author}{T.~Dupont},
\newblock \bibinfo{title}{Thin metamaterial using acoustic black hole profiles for broadband sound absorption},
\newblock \bibinfo{journal}{Applied Acoustics} \bibinfo{volume}{216} (\bibinfo{year}{2024}) \bibinfo{pages}{109744}.
\bibitem[{Zhang et~al.(2025)Zhang, Geng, Zhao, Cao, and Wang}]{zhang2025multi}
\bibinfo{author}{X.~Zhang}, \bibinfo{author}{M.~Geng}, \bibinfo{author}{C.~Zhao}, \bibinfo{author}{Y.~Cao}, \bibinfo{author}{P.~Wang},
\newblock \bibinfo{title}{Multi-gradient acoustic black hole metamaterial for near-perfect sound attenuation: Theory, simulation and experiments},
\newblock \bibinfo{journal}{Applied Acoustics} \bibinfo{volume}{231} (\bibinfo{year}{2025}) \bibinfo{pages}{110546}.
\bibitem[{Foucaud et~al.(2014)Foucaud, Michon, Gourinat, Pelat, and Gautier}]{foucaud2014artificial}
\bibinfo{author}{S.~Foucaud}, \bibinfo{author}{G.~Michon}, \bibinfo{author}{Y.~Gourinat}, \bibinfo{author}{A.~Pelat}, \bibinfo{author}{F.~Gautier},
\newblock \bibinfo{title}{Artificial cochlea and acoustic black hole travelling waves observation: Model and experimental results},
\newblock \bibinfo{journal}{Journal of Sound and Vibration} \bibinfo{volume}{333} (\bibinfo{year}{2014}) \bibinfo{pages}{3428--3439}.
\bibitem[{Rao(2006)}]{RAO}
\bibinfo{author}{S.~Rao}, \bibinfo{title}{Approximate Analytical Methods}, \bibinfo{publisher}{John Wiley \& Sons, Ltd}, \bibinfo{year}{2006}, pp. \bibinfo{pages}{647--699}. \URLprefix \url{https://onlinelibrary.wiley.com/doi/abs/10.1002/9780470117866.ch17}. \DOIprefix\doi{https://doi.org/10.1002/9780470117866.ch17}. \href{http://arxiv.org/abs/https://onlinelibrary.wiley.com/doi/pdf/10.1002/9780470117866.ch17}{\tt arXiv:https://onlinelibrary.wiley.com/doi/pdf/10.1002/9780470117866.ch17}.
\bibitem[{E756--05(2023)}]{E756}
E756--05, \bibinfo{title}{Standard Test Method for Measuring Vibration-Damping Properties of Materials}, \bibinfo{type}{Standard}, ASTM International, \bibinfo{address}{West Conshohocken, PA}, \bibinfo{year}{2023}.
\bibitem[{Brun et~al.(2022)Brun, Cortés, García-Barruetabeña, Sarría, and Elejabarrieta}]{Brun_2022}
\bibinfo{author}{M.~Brun}, \bibinfo{author}{F.~Cortés}, \bibinfo{author}{J.~García-Barruetabeña}, \bibinfo{author}{I.~Sarría}, \bibinfo{author}{M.~J. Elejabarrieta},
\newblock \bibinfo{title}{A robust technique for polymer damping identification using experimental transmissibility data},
\newblock \bibinfo{journal}{Polymers} \bibinfo{volume}{14} (\bibinfo{year}{2022}) \bibinfo{pages}{2535}.
\bibitem[{Soroor et~al.(2021)Soroor, Asgari, and Haddadpour}]{soroor2021effect}
\bibinfo{author}{A.~O. Soroor}, \bibinfo{author}{M.~Asgari}, \bibinfo{author}{H.~Haddadpour},
\newblock \bibinfo{title}{Effect of axially graded constraining layer on the free vibration properties of three layered sandwich beams with magnetorheological fluid core},
\newblock \bibinfo{journal}{Composite Structures} \bibinfo{volume}{255} (\bibinfo{year}{2021}) \bibinfo{pages}{112899}.
\bibitem[{Denis et~al.(2015)Denis, Gautier, Pelat, and Poittevin}]{denis2015measurement}
\bibinfo{author}{V.~Denis}, \bibinfo{author}{F.~Gautier}, \bibinfo{author}{A.~Pelat}, \bibinfo{author}{J.~Poittevin},
\newblock \bibinfo{title}{Measurement and modelling of the reflection coefficient of an acoustic black hole termination},
\newblock \bibinfo{journal}{Journal of Sound and Vibration} \bibinfo{volume}{349} (\bibinfo{year}{2015}) \bibinfo{pages}{67--79}.
\bibitem[{Hook et~al.(2019)Hook, Cheer, and Daley}]{hook2019}
\bibinfo{author}{K.~Hook}, \bibinfo{author}{J.~Cheer}, \bibinfo{author}{S.~Daley},
\newblock \bibinfo{title}{A parametric study of an acoustic black hole on a beam},
\newblock \bibinfo{journal}{The Journal of the Acoustical Society of America} \bibinfo{volume}{145} (\bibinfo{year}{2019}) \bibinfo{pages}{3488--3498}.
\bibitem[{Gustavsen and Semlyen(1999)}]{gustavsen_rational_1999}
\bibinfo{author}{B.~Gustavsen}, \bibinfo{author}{A.~Semlyen},
\newblock \bibinfo{title}{Rational approximation of frequency domain responses by vector fitting},
\newblock \bibinfo{journal}{IEEE Transactions on Power Delivery} \bibinfo{volume}{14} (\bibinfo{year}{1999}) \bibinfo{pages}{1052--1061}.
\bibitem[{Gustavsen(2006)}]{gustavsen_improving_2006}
\bibinfo{author}{B.~Gustavsen},
\newblock \bibinfo{title}{Improving the pole relocating properties of vector fitting},
\newblock \bibinfo{journal}{IEEE Transactions on Power Delivery} \bibinfo{volume}{21} (\bibinfo{year}{2006}) \bibinfo{pages}{1587--1592}.
\bibitem[{Deschrijver et~al.(2008)Deschrijver, Mrozowski, Dhaene, and Zutter}]{deschrijver_macromodeling_2008}
\bibinfo{author}{D.~Deschrijver}, \bibinfo{author}{M.~Mrozowski}, \bibinfo{author}{T.~Dhaene}, \bibinfo{author}{D.~D. Zutter},
\newblock \bibinfo{title}{Macromodeling of {Multiport} {Systems} {Using} a {Fast} {Implementation} of the {Vector} {Fitting} {Method}},
\newblock \bibinfo{journal}{IEEE Microwave and Wireless Components Letters} \bibinfo{volume}{18} (\bibinfo{year}{2008}) \bibinfo{pages}{383--385}.
\bibitem[{Sangle et~al.(2024)Sangle, Alajlouni, and Tarazaga}]{sangle2024pole}
\bibinfo{author}{S.~Sangle}, \bibinfo{author}{S.~Alajlouni}, \bibinfo{author}{P.~A. Tarazaga},
\newblock \bibinfo{title}{A pole-based approach to interpret electromechanical impedance measurements in structural health monitoring},
\newblock \bibinfo{journal}{arXiv preprint arXiv:2411.05871}  (\bibinfo{year}{2024}).
\bibitem[{Austin and Cheer(2022)}]{austin2022realisation}
\bibinfo{author}{B.~Austin}, \bibinfo{author}{J.~Cheer},
\newblock \bibinfo{title}{Realisation of acoustic black holes using multi-material additive manufacturing},
\newblock \bibinfo{journal}{Frontiers in Physics} \bibinfo{volume}{10} (\bibinfo{year}{2022}) \bibinfo{pages}{1070345}.
\bibitem[{Yang et~al.(2025)Yang, Yu, and Ye}]{yang2025study}
\bibinfo{author}{C.~Yang}, \bibinfo{author}{H.~Yu}, \bibinfo{author}{T.~Ye},
\newblock \bibinfo{title}{Study on broadband vibration reduction characteristics and optimal design of the acoustic black hole plate with damping oscillators},
\newblock \bibinfo{journal}{Scientific Reports} \bibinfo{volume}{15} (\bibinfo{year}{2025}) \bibinfo{pages}{12100}.

\end{thebibliography}


\bio{}
\endbio

\end{document}